%15.1.08
%math.GT/0603429  i6daj
\input amstex
\documentstyle{amsppt}
\NoBlackBoxes
\magnification=\magstep1
\advance\vsize-0.5cm\voffset=-1cm\advance\hsize1cm\hoffset0cm
\define\R{{\Bbb R}} \define\Z{{\Bbb Z}} \define\Q{{\Bbb Q}} \define\C{{\Bbb C}}
\def\pr{\mathop{\fam0 pr}}
\def\forg{\mathop{\fam0 forg}}

\def\Sq{\mathop{\fam0 Sq}}

\def\Emb{\mathop{\fam0 Emb}}

\def\id{\mathop{\fam0 id}}

\def\link{\mathop{\fam0 link}}

\def\Int{\mathop{\fam0 Int}}

\def\Cl{\mathop{\fam0 Cl}}
\def\Emb{\mathop{\fam0 Emb}}

\def\im{\mathop{\fam0 im}}

\def\t{\widetilde}

\newpage
\topmatter
\title A classification of smooth embeddings of 3-manifolds in 6-space 
\endtitle
\author Arkadiy Skopenkov \endauthor
\address
Department of Differential Geometry, Faculty of Mechanics and
Mathematics, Moscow State University, 119992, Moscow, Russia,
and Independent University of Moscow, B. Vlasyevskiy, 11, 119002,
Moscow, Russia.
e-mail: skopenko\@mccme.ru \endaddress
\thanks
The original publication is available at www.springerlink.com
DOI 10.1007/s00209-007-0294-1. 
\newline
The author gratefully acknowledges the support from Deligne 2004 Balzan
prize in mathematics, the Russian Foundation for Basic Research Grants 
07-01-00648a and 06-01-72551-NCNILa, President of Russian Federation Grants 
MD-4729.2007.1 and NSH-4578.2006.1, 
DAAD Stipendium Leonhard Euler Projekt 05/05857.
\endthanks

\abstract
We work in the smooth category.
If there are knotted embeddings $S^n\to\R^m$, which often happens for
$2m<3n+4$, then no explicit complete description of the embeddings of 
$n$-manifolds into $\R^m$ up to isotopy was known, except for the disjoint 
unions of spheres.
Let $N$ be a closed connected orientable 3-manifold.
Our main result is the following description of the set $\Emb^6(N)$ of
embeddings $N\to\R^6$ up to isotopy.
We define the Whitney and the Kreck invariants and prove that 
{\it the Whitney invariant $W:\Emb^6(N)\to H_1(N;\Z)$ is surjective.
For each $u\in H_1(N;\Z)$ the Kreck invariant $\eta_u:W^{-1}u\to\Z_{d(u)}$
is bijective, where $d(u)$ is the divisibility of the projection of $u$ to 
the free part of $H_1(N;\Z)$.}

The group $\Emb^6(S^3)$ is isomorphic to $\Z$ (Haefliger). 
This group acts on $\Emb^6(N)$ by embedded connected sum. 
It was proved that the orbit space of this action maps under $W$ bijectively to 
$H_1(N;\Z)$ (by Vrabec and Haefliger's smoothing theory). 
The new part of our classification result is {\it the determination of the 
orbits} of the action.
E. g. for $N=\R P^3$ the action is free, while for $N=S^1\times S^2$ we
explicitly construct {\it an embedding $f:N\to\R^6$ such that for each knot
$l:S^3\to\R^6$ the embedding $f\# l$ is isotopic to $f$}.

The proof uses new approaches involving modified surgery theory as developed by 
Kreck or the Bo\'echat-Haefliger formula for the smoothing obstruction.
\endabstract
\endtopmatter

\document
\head 1. Introduction \endhead

{\bf Main results. }

This paper is on the classical {\it Knotting Problem} in topology: 
{\it given an $n$-manifold $N$ and a number $m$, describe
isotopy classes of embeddings $N\to\R^m$}.
For recent surveys see [RS99, Sk07].
We work in the smooth category. 
Let $\Emb^m(N)$ be the set of embeddings $N\to\R^m$ up to isotopy. 
% (or, which is the same, up to orientation-preserving isoposition).

The Knotting Problem is more accessible for $2m\ge3n+4$.  
%For $2m\ge3n+4$ the set $\Emb^m(N)$ was described using the {\it Haefliger-Wu
%invariant} and the invariants that could be derived from it. For many cases 
%this description allows explicit calculations [Ha63, RS99, \S4, Sk02, \S1, Sk07, \S5].
The Knotting Problem is much harder for 
$$2m\le3n+3:$$
if $N$ is a closed manifold that is not a disjoint union of spheres, then no 
explicit complete descriptions of isotopy classes was known
\footnote{See though [KS05, Sk06]; for {\it rational} and {\it piecewise 
linear} classification see [CRS07', CRS] and [Hu69, \S12, Vr77, Sk97, Sk02, 
Sk05, Sk06], respectively.}, in spite of the 
existence of interesting approaches [Br68, Wa70, GW99]\footnote{
I am grateful to M. Weiss for indicating that the approach of [GW99, We] does 
give explicit results on higher homotopy groups of the space of embeddings 
$S^1\to\R^n$.}. 

In particular, the classification of embeddings $N\to\R^{2n}$ was known for
$n\ge4$ [HH63, Ba75, Vr77, Sk07, \S2] and was unknown for $n=3$ (except
for a disjoint union of 3-spheres).

In this paper we address the case $(m,n)=(6,3)$ and, more
generally, $2m=3n+3$. We assume everywhere that

{\it $N$ is a closed connected orientable 3-manifold,}

unless the contrary is explicitly indicated.
Our main result is a complete explicit description of the set $\Emb^6(N)$ of
embeddings $N\to\R^6$ up to isotopy.

We omit $\Z$-coefficients from the notation of (co)ho\-mo\-lo\-gy
groups.
We define the Whitney and the Kreck invariants below in \S1.

\smallskip
{\bf Classification Theorem.} {\it For every closed connected
orientable 3-manifold $N$ the Whitney invariant
$$W:\Emb\phantom{}^6(N)\to H_1(N)$$
is surjective.
For each $u\in H_1(N)$ the Kreck invariant
$$\eta_u:W^{-1}(u)\to\Z_{d(u)}$$
is  bijective, where $d(u)$ is the divisibility of the projection of $u$ to
the free part of $H_1(N)$.}

\smallskip
Recall that for an abelian group $G$ the divisibility of zero is zero and the 
divisibility of $x\in G-\{0\}$ is 
%\linebreak 
$\max\{d\in\Z\ | \ \text{there is }x_1\in G:\ x=dx_1\}$.

\smallskip
{\bf Corollary.} {\it (1) The Kreck invariant $\eta_0:\Emb^6(N)\to\Z$ is a 1--1 
correspondence if $N=S^3$ [Ha66] or an integral homology sphere [Ha72, Ta06]. 

(2) If $H_2(N)=0$ (i.e. $N$ is a rational homology sphere, e.g. $N=\R P^3$), 
then $\Emb^6(N)$ is in (non-canonical) 1--1 correspondence with 
$\Z\times H_1(N)$. 

(3) Embeddings $S^1\times S^2\to\R^6$ with zero Whitney invariant are in 
1--1 correspondence with $\Z$, and for each integer $k\ne0$ there are exactly 
$k$ distinct embeddings $S^1\times S^2\to\R^6$ with the Whitney invariant $k$, 
cf. Corollary (a) and (b) below. 

(4) The Whitney invariant 
$W:\Emb^6(N_1\# N_2)\to H_1(N_1\# N_2)\cong H_1(N_1)\oplus H_1(N_2)$ is 
surjective and $\# W^{-1}(a_1\oplus a_2)=GCD(\# W_1^{-1}(a_1),\# W_2^{-1}(a_2))$, 
where $GCD(\infty,a)=a$.}   
  
\smallskip
All isotopy classes of embedings $N\to\R^6$ can be simply constructed (from 
a certain given embedding), see the end of \S1. 

Note that some 3-manifolds appear in the theory of integrable systems together 
with their embeddings into $\R^6$ (given by a system of algebraic equations) 
[BF04, Chapter 14]. For other examples of embeddings see %[Sk06, 
beginning of \S5. 
%].  

Notice that our classification of smooth embeddings $N\to\R^6$ is similar to 
the Wu classification of smooth immersions $N\to\R^5$ [SST02, Theorem 3.1] and 
to the Pontryagin classification of 
homotopy classes of maps $N\to S^2$ (or, which is the same, non-zero vector 
fields on $N$) [CRS07]. 
It would be interesting to construct maps $[N;S^2]\leftrightarrow\Emb^6(N)$. 

By Corollary (1) our achievement is the transition from $N=S^3$ to arbitrary 
(closed connected orientable) 3-manifolds.
Let us explain what is involved. 
It was known that the embedded connected sum defines a group structure on 
$\Emb^6(S^3)$ and that $\Emb^6(S^3)\cong\Z$ [Ha66]. 
The group $\Emb^6(S^3)$ acts on the set 
$\Emb^6(N)$ by connected summation of embeddings $g:S^3\to\R^6$ and 
$f:N\to\R^6$ whose images are contained in disjoint cubes. 
By general position the connected sum $f\#g$ is well-defined, i.e. does not 
depend on the choice of an arc between $gS^3$ and $fN$.
It was  known that the orbit space of this action maps under $W$ (defined in
a different way) bijectively to $H_1(N)$. 
\footnote{This follows from the PL Classification Theorem and the definition of 
the Kreck invariant in \S4, cf. [Ta06, Proposition 2.4].}
Thus the knotting problem was reduced to {\it the determination of the orbits
of this action}, which is as non-trivial a problem and the new part of 
the Classification Theorem.

%Denote $kt=t\#\dots\#t$ ($k$ summands). 

\smallskip
{\bf Addendum to the Classification Theorem.}
{\it If $f:N\to\R^6$ is an embedding, $t$ is the generator of
$\Emb^6(S^3)\cong\Z$ and $kt$ is a connected sum of $k$ copies of $t$, then
$$W(f\#kt)=W(f)\quad\text{and}\quad
\eta_{W(f)}(f\#kt)\equiv\eta_{W(f)}(f)+k\mod d(W(f)).$$}
Here the first equality follows by the definition of the Whitney invariant 
(see below), and the second equality is proved in the subsection 
`definition of the Kreck invariant'.

E. g. for $N=\R P^3$ the action of $\Emb^6(S^3)$ on $\Emb^6(N)$ is free 
(because $H_1(\R P^3)\cong\Z_2$) while for $N=S^1\times S^2$ we have the 
following corollary.
\footnote{We believe that this very concrete corollary or the case $N=\R P^3$
are as non-trivial as the general case of the Classification Theorem.}

\smallskip
{\bf Corollary.}
{\it (a) There is an embedding $f:S^1\times S^2\to\R^6$ such that for
each knot $l:S^3\to\R^6$ the embedding $f\# l$ is isotopic to $f$.

Take the standard embeddings $2D^4\times S^2\subset\R^6$ (where 2 is 
multiplication by 2) and $\partial D^2\subset\partial D^4$. 
Fix  a point $x\in S^2$. 
Such an embedding $f$ is the connected sum of 
$$2\partial D^4\times x\quad\text{with}\quad
\partial D^2\times S^2\subset D^4\times S^2\subset 2D^4\times S^2\subset\R^6.$$ 
(b) For each embedding $f:N\to\R^6$ such that $f(N)\subset\R^5$ (e.g. for the 
standard embedding $f:S^1\times S^2\to\R^6$) and each non-trivial knot 
$l:S^3\to\R^6$ the embedding $f\# l$ is not isotopic to $f$.}

\smallskip
Corollary (a) is generalized and proved at the end of \S1.
Corollary (b) follows from the Classification Theorem and (the easy necessity 
of $2W(f)=0$ of) the Compression Theorem stated below in \S1.

In \S2 and \S4 we present two proofs of the Classification Theorem. 
These sections are independent of each other, except for the common reduction 
of the Classification Theorem (to the Injectivity Lemma) and the use in \S4 of 
$\eta:\Emb^6(S^3)\to\Z$ defined at the beginning of \S2.  
In order to let the reader understand the main ideas before going into details, 
we sometimes apply a lemma before presenting its proof. 
 
The two proofs correspond to the two definitions of the Kreck invariant (given 
in \S1 and in \S4). 
We present both proofs because their ideas and generalizations are distinct.

The first proof (\S2) is the result of a discussion with Matthias Kreck who 
kindly allowed it to be included it in this paper.
It uses a new approach involving modified surgery as developed by Kreck [Kr99]; 
for ideas of this proof see [KS05, \S1].
This proof is self-contained in the sense that it does not use results outside 
modified surgery.
In particular we reprove the injectivity of the Haefliger invariant 
$\Emb^6(S^3)\to\Z$. 
\footnote{The surjectivity can be reproven analogously, cf. [Fu94]; 
for an alternative reproof of the surjectivity see [Ta04].}
This approach is useful in other relatively low dimensions [KS05].

The second proof (\S4) uses many results: the Haefliger construction of the 
isomorphism $\Emb^6(S^3)\cong\Z$, Haefliger smoothing theory [Ha67] and 
the Bo\'echat-Haefliger smoothing result [BH70] (the latter uses either the 
calculation of cobordism classes of PL embeddings of 4-manifolds into $\R^7$ 
[BH70] or the argument of [Bo71, p. 153]).

The second proof generalizes to the following result proved in \S4, cf. 
[Sk07, \S2 and \S3].

\smallskip
{\bf Higher-dimensional Classification Theorem.}

{\it (a) $\Emb\phantom{}^{6k}(S^p\times S^{4k-1-p})\cong
\pi_{4k-1-p}(V_{2k+p+1,p+1})\oplus\Z$ for $1\le p\le2k-2$, where $V_{a,b}$ is 
the Stiefel manifold of $b$-frames in $\R^a$.

(b) Let $N$ be a closed homologically $(2k-2)$-connected $(4k-1)$-manifold.  
%such that $w_{2k}(N)=0$.
Then the Whitney invariant $W:\Emb^{6k}(N)\to H_{2k-1}(N)$ is surjective
and for each $u\in H_{2k-1}(N)$ there is a bijective invariant
$\eta_u:W^{-1}u\to\Z_{d(u)}$.}

\smallskip
The following particular case of (b) should be compared with (a): 

{\it The Whitney invariant
$W:\Emb\phantom{}^{6k}(S^{2k-1}\times S^{2k})\to\Z$ is
surjective and for each $u\in\Z$ there is a bijective invariant
$\eta_u:W^{-1}u\to\Z_u$.} 

An alternative proof of the Higher-dimensional Classification Theorem (a) for 
$k>1$ can be obtained using the construction of the smooth Whitehead torus 
[Sk06]; such an argument is also non-trivial.

\bigskip
{\bf Definition of the Whitney invariant.}
\nopagebreak

Fix orientations on $\R^6$ and on $N$. 
Let $B^3$ be a closed 3-ball in $N$.
Denote $N_0:=\Cl(N-B^3)$.
From now on $f:N\to\R^6$ is an embedding, unless another meaning of $f$ 
is explicitly given.   

Since $N$ is orientable, $N$ embeds into $\R^5$ [Hi61].
Fix an embedding $g:N\to\R^6$ such that $g(N)\subset\R^5$.
The restrictions of $f$ and $g$ to $N_0$ are regular homotopic by [Hi60].
Since $N_0$ has a 2-dimensional spine, it follows that these restrictions
are isotopic, cf. [HH63, 3.1.b, Ta06, Lemma 2.2].
So we can make an isotopy of $f$ and assume that $f=g$ on $N_0$.
Take a general position homotopy $F:B^3\times I\to\R^6$ relative to
$\partial B^3$ between the restrictions of $f$ and $g$ to $B^3$.
Then $f\cap F:=(f|_{N-B^3})^{-1}F(B^3\times I)$ (i.e. `the intersection of this
homotopy with $f(N-B^3)$') is a 1-manifold (possibly non-compact) without 
boundary.
Define $W(f)$ to be the homology class of the closure of this oriented 
1-manifold:
$$W(f):=[\Cl(f\cap F)]\in H_1(N_0,\partial N_0)\cong H_1(N).$$
The orientation on $f\cap F$ is defined as follows. 
For each point $x\in f\cap F$ take a vector at $x$ tangent to
$f\cap F$.
Complete this vector to a positive base tangent to $N$.
By general position there is a unique point $y\in B^3\times I$ such that 
$Fy=fx$.
The tangent vector at $x$ thus gives a tangent vector at $y$ to $B^3\times I$.
Complete this vector to a positive base tangent to $B^3\times I$, where the 
orientation on $B^3$ comes from $N$.
The union of the images of the constructed two bases is a base at $Fy=fx$ 
of $\R^6$.
If this base is positive, then call the initial vector of $f\cap F$ positive.
Since a change of the orientation on $f\cap F$ forces a change of the 
orientation of the latter base of $\R^6$, it follows that this condition indeed 
defines an orientation on $f\cap F$.

\smallskip
The Whitney invariant is well-defined, i.e. independent of the choice of
$F$ and of the isotopy making $f=g$ outside $B^3$.
This is so because the above definition is clearly equivalent to that of
[Hu69, \S12, Vr77, p. 145, Sk07, \S2] (analogous to the definition of the
Whitney obstruction to embeddability): $W(f)$ is the homology class of the
algebraic sum of the top-dimensional simplices of the self-intersection set
$\Sigma(H):=\Cl\{x\in N\times I\ |\ \#H^{-1}Hx>1\}$ of a general position
homotopy $H$ between $f$ and $g$ (for definition of the signs of the simplices
see [Hu69, \S12, Vr77, p. 145, Sk07, \S2]).
For another equivalent definition see [HH63, bottom of p. 130 and p. 134],
cf. the Bo\'echat-Haefliger Invariant Lemma of \S2.

Clearly, $W(g)=0$.
(Corollary (b) for $f=g$ follows from the Classification Theorem and $W(g)=0$.)  

The definition of $W$ depends on the choice of $g$, but we write $W$ not
$W_g$ for brevity.
If $H_1(N)$ has no 2-torsion, then $W$ is in fact independent of the choice of 
$g$ by (the necessity of $2W(f)=0$ of) the Compression Theorem at the end of 
\S1 or by the Bo\'echat-Haefliger Invariant Lemma of \S2.

Since a change of the orientation on $N$ forces a change of the orientation on 
$B^3$, the class $W(f)$ is independent of the choice of the orientation on $N$.
For the reflection $\sigma:\R^6\to\R^6$ with respect to a hyperplane 
we have $W(\sigma\circ f)=-W(f)$ 
(because we may assume that $f=g=\sigma\circ f$ on $N_0$ and because a change 
of the orientation of $\R^6$ forces a change of the orientation of $f\cap F$). 

For a closed homologically $(2k-2)$-connected $(4k-1)$-manifold $N$ 
the Whitney invariant
$W:\Emb^{6k}(N)\to H_{2k-1}(N)$ is defined analogously to the above by
$$W(f):=[\Cl(f|_{N-B^{4k-1}})^{-1}F(B^{4k-1}\times I)]
\in H_{2k-1}(N_0,\partial N_0)\cong H_{2k-1}(N).$$

\smallskip
%\bigskip
{\bf Definition of the Kreck invariant.}
\nopagebreak

Denote by 

$\bullet$
$C_f$ the closure of the complement in $S^6\supset\R^6$ to a tubular
neighborhood of $f(N)$ and 

$\bullet$
$\nu_f:\partial C_f\to N$ the restriction of the normal bundle of $f$.  

An orientation-preserving diffeomorphism 
$\varphi:\partial C_f\to\partial C_{f'}$ such that $\nu_f=\nu_{f'}\varphi$ 
is simply called an {\it isomorphism}. 
For an isomorphism $\varphi$ denote 
$$M=M_\varphi:=C_f\cup_\varphi(-C_{f'}).$$ 
An isomorphism $\varphi:\partial C_f\to\partial C_{f'}$ is called {\it spin}, 
if $\varphi$  over $N_0$ is defined by an isotopy between the restrictions of 
$f$ and $f'$ to $N_0$.  
A spin isomorphism exists because the restrictions to $N_0$ of $f$ and $f'$ are 
isotopic (see the definition of the Whitney invariant) and $\pi_2(SO_3)=0$. 

\smallskip
{\bf Spin Lemma.} 
{\it If $\varphi$ is a spin isomorphism, then $M_\varphi$ is spin.} 

\smallskip
{\it Proof.}
By a {\it spin structure}  we mean a {\it stable spin structure}.
Recall that a (stable) spin structure on $M$ is a framing of the stable 
normal bundle on the 2-skeleton of $M$, and that stable spin structures
are equivalent if they are equivalent over the 1-skeleton of $M$ [Ki89, IV].
Take spin structures on $C_f$ and $C_{f'}$ obtained from their embeddings into 
$S^6$ and standard normal framings $S^6\subset S^{13}$. 

The 2-skeleton of $\partial C_f\cong N\times S^2$ is contained in
$\nu_f^{-1}N_0$.
Consider the spin structure on $\partial C_f$ induced by the {\it framed} 
embedding $\partial C_f\subset S^6$ (the framing is the normal vector field 
looking to the connected component of $S^6-\partial C_f$ containing $N$). 
Consider an analogous spin structure on $\partial C_{f'}$. 
By the definition of $\varphi|_{N_0}$ the first spin structure goes to the 
second one under $\varphi|_{N_0}$. 
Hence the spin structures on $C_f$ and $C_{f'}$ agree on the boundary. 
So the manifold $M$ is spin. 
\qed

\smallskip
Denote by $\sigma (X)$ the signature of a 4-manifold $X$.  
Denote by $PD:H^i(Q)\to H_{q-i}(Q,\partial Q)$ and 
$PD:H_i(Q)\to H^{q-i}(Q,\partial Q)$ Poincar\'e duality (in any manifold $Q$). 
For $y\in H_4(M_\varphi)$ and a $k$-submanifold $C\subset M_\varphi$ (e.g. $C=C_f$ or 
$C=\partial C_f$) denote 
$$y\cap C:=PD[(PDy)|_C]\in H_{k-2}(C,\partial C).$$ 
If $y$ is represented by a closed oriented 4-submanifold $Y\subset M_\varphi$ in 
general position to $C$, then $y\cap C$ is represented by $Y\cap C$.

A {\it homology Seifert surface} for $f$ is the image $A_f$ of the fundamental 
class $[N]$ under the composition $H_3(N)\to H^2(C_f)\to H_4(C_f,\partial C_f)$ 
of the Alexander and Poincar\'e duality isomorphisms, cf. \S3. 
\footnote{This seems to be the only notion of this subsection whose analogue for knots 
$S^1\subset S^3$ is useful.}

A {\it joint homology Seifert surface} for $f$ and $f'$ is a class 
$y\in H_4(M_\varphi)$ such that 
$$y\cap C_f=A_f\quad\text{and}\quad y\cap C_{f'}=A_{f'}.$$  

%\smallskip
{\bf Agreement Lemma.} 
{\it If $\varphi$ is a spin isomorphism and $W(f)=W(f')$, then there is a joint 
homology Seifert surface for $f$ and $f'$.}

\smallskip
The proof (\S2) is non-trivial and is an essential step in the proof of the 
the Classification Theorem. 
\footnote{Although we do not need it, we note that if $\varphi$ is a spin isomorphism 
and there is a joint homology Seifert surface for $f$ and $f'$, then 
$W(f)=W(f')$.}  

We identify with $\Z$ the zero-dimensional homology groups and the 
$n$-dimensional cohomology groups of closed connected oriented $n$-manifolds. 
The intersection products in 6-manifolds are omitted from the notation. 
For a closed connected oriented 6-manifold $Q$ and $x\in H_4(Q)$ denote by 
$$\sigma_x(Q):=\frac{xPDp_1(Q)-x^3}3\in H_0(Q)=\Z$$ 
the virtual signature of $(Q,x)$. 
(Since $H_4(Q)\cong[Q,\C P^\infty]$, there is a closed connected 
oriented 4-submanifold $X\subset Q$ representing the class $x$. 
Then by [Hi66', end of 9.2] or else by the Submanifold Lemma below we have 
$3\sigma(X)=p_1(X)=xPDp_1(Q)-x^3=3\sigma_x(Q)$.) 

{\it The Kreck invariant} of two embeddings $f$ and $f'$ such that $W(f)=W(f')$ 
is defined by  
$$\eta(f,f'):\equiv\frac{\sigma_{2y}(M_\varphi)}{16}=
\frac{yPDp_1(M_\varphi)-4y^3}{24}\mod d(W(f)),$$
where $\varphi:\partial C_f\to\partial C_{f'}$ is a spin isomorphism and 
$y\in H_4(M)$ is a joint homology Seifert surface for $f$ and $f'$. 
Cf. [Ek01, 4.1, Zh]. 
We have $2y\mod2=0=PDw_2(M)$, so any closed connected oriented 4-submanifold of 
$M$ representing the class $2y$ is spin, hence by the Rokhlin Theorem 
$\sigma_{2y}(M)$ is indeed divisible by 16. 

The Kreck invariant is well-defined by the following 

\smallskip
{\bf Independence Lemma.}
{\it The residue $\sigma_{2y}(M_\varphi)/16\mod d(W(f))$ is independent of the 
choice of the spin isomorphism $\varphi$ and of the joint homology Seifert 
surface $y$.}

\smallskip
The proof (\S2) is non-trivial and is an essential step in the proof of the 
Classification Theorem. 

For $u\in H_1(N)$ fix an embedding $f':N\to\R^6$ such that $W(f')=u$ and define 
$\eta_u(f):=\eta(f,f')$. 
(We write $\eta_u(f)$ not $\eta_{f'}(f)$ for simplicity.) 

The choice of the other orientation for $N$ (resp. $\R^6$) will in general give 
rise to different values for the Kreck invariant. 
But such a choice only permutes the bijection $W^{-1}(u)\to\Z_{d(u)}$ (resp. 
replaces it with the bijection $W^{-1}(-u)\to\Z_{d(u)}$). 

\smallskip
{\it Proof of the second equality of the Addendum.} 
We may assume that $\nu_f=\nu_{f\#t}$ outside $B^3$. 
Take a spin isomorphism $\varphi:\partial C_f\to\partial C_{f\#t}$ that is 
identical outside $B^3$. 

We have $C_{f\#t}=C_f\natural C_t$, where the last boundary connected sum is 
along an isomorphism $\nu_f^{-1}B^3\to \nu_t^{-1}B^3$. 
A proper connected orientable 4-submanifold $X\subset M_{\id}$ representing 
twice the joint homology Seifert surfaces for $f$ and $f$ (which need not be 
empty) intersects $\nu_f^{-1}B^3$ by $B^3\times*$, where $*\in S^2$. 
We can take a proper connected orientable 4-submanifold 
$X_t\subset M$ representing twice the joint homology Seifert surfaces for $t$ 
and the standard embedding $S^3\to S^6$, that intersects $\nu_t^{-1}B^3$ by 
$B^3\times*$, where $*\in S^2$. 

Hence $X\#X_t\subset M_\varphi$ represents twice the joint homology Seifert 
surfaces for $f$ and $f\#t$. 
Now the second equality of the Addendum follows because 
$$16\eta(f\#t,f)\underset{\mod16d(W(f))}\to\equiv
\sigma(X\#X_t)=\sigma(X)+\sigma(X_t)=16\eta_0(t,g_0)=16,$$
where $g_0:S^3\to\R^6$ is the standard embedding. 
In this formula the second equality holds by the Novikov-Rokhlin additivity. 
The fourth equality holds by the Kreck Invariant Lemma below and 
[GM86, Remarks to the four articles of Rokhlin, II.2.7 and III.excercises.IV.3, 
Ta04] (or by 
[Ha66] because the Kreck invariant for $N=S^3$ coincides with the Haefliger 
invariant, see the very beginning of \S2, cf. %[Sk06, 
\S5, (4)). 
\qed

\smallskip
Let us present a formula for the Kreck invariant analogous to 
[GM86, Remarks to the four articles of Rokhlin, II.2.7 and III.excercises.IV.3, 
Ta04, Corollary 6.5, Ta06, Proposition 4.1]. 
This formula is useful when an embedding goes through $\R^5$ (cf. the next 
subsection) or is given by a system of equations (because we can obtain a 
`Seifert surface' by changing the equality to the inequality in one 
of the equations). 

\smallskip
{\bf The Kreck Invariant Lemma.} {\it Let $f,f':N\to\R^6$ be two embeddings 
such that $W(f)=W(f')$, $\varphi:\partial C_f\to\partial C_{f'}$ a spin 
isomorphism, 
$Y\subset M_\varphi$ a closed connected oriented 4-submanifold representing a 
joint homology Seifert surface and $\overline p_1\in\Z$, 
$\overline e\in H_2(Y)$ are the Pontryagin number and Poincar\'e dual of the 
Euler classes of the normal bundle of $Y$ in $M_\varphi$. 
Then  
$$\frac{\sigma_{2[Y]}(M_\varphi)}{16}=\frac{\sigma(Y)-\overline p_1}8=
\frac{\sigma(Y)-\overline e\cap\overline e}8.$$}

%\smallskip
{\bf Submanifold Lemma.} {\it Let $Y$ be a closed oriented connected 
4-submanifold of a closed orientable connected 6-manifold $Q$. 
Denote by $[Y]\in H_4(Q)$ the homology class of $Y$, by 
$\overline p_1\in\Z$ and $\overline e\in H_2(Y)$ the Pontryagin number and 
Poincar\'e dual of the Euler class of the normal bundle of $Y$ in $Q$. 
Then 
$$[Y]^3=\overline e\cap\overline e=\overline p_1\quad\text{and}
\quad [Y]PDp_1(Q)=p_1(Y)+\overline p_1.$$}

%\smallskip
{\it Proof (folklore).} 
We have $[Y]PDp_1(Q)=p_1(Q)|_Y=p_1(Y)+\overline p_1$, where the last  
equality holds because $\tau_Q|_Y\cong\tau_Y\oplus\nu_Q(Y)$. 

%Let $Y_1,Y_2\subset Q$ be submanifolds close to $Y$ and in general position 
%to each other and to $Y$. Then
We have 
$$[Y]^3=([Y]\cap Y)\bigcap\limits_Y([Y]\cap Y)=
\overline e\bigcap\limits_Y\overline e=\overline p_1.$$
%Here $\#$ means the algebraic sum of intersection points in $Q$, and 
%$Y_1\cap Y$, $Y_2\cap Y$ mean {\it oriented} intersections.  
In order to prove the latter equality take two general position sections of 
the normal bundle of $Y\subset Q$. 
Then $\overline e$ is the homology class of the 
appropriately oriented zero submanifold of any section. 
The class $\overline p_1$ is the homology class of the appropriately oriented 
submanifold on which the rank of pair of normal vectors formed by the two 
sections is less than 1, i.e. on which {\it both} vectors are zeros. 
Thus $\overline e\cap\overline  e=\overline p_1$.  
\qed

\smallskip
The Kreck Invariant Lemma holds by the Submanifold Lemma  
because 
$$3\sigma_{2y}(M_\varphi)=2yPDp_1(M_\varphi)-8y^3=
2p_1(Y)+2\overline p_1-8\overline p_1=6\sigma(Y)-6\overline p_1.$$ 
We remark that the assumption $W(f)=W(f')$ of the Kreck Invariant Lemma follows 
from the other assumptions, but we anyway use the Lemma in situations when we 
know that $W(f)=W(f')$. 

%For a closed homologically $(2k-2)$-connected $(4k-1)$-manifold $N$ the Kreck 
%invariant is defined analogously to the above. No! X not embedded
%See an equivalent definition of the Kreck invariant in \S4. 

\bigskip
{\bf The Compression Problem. }

{\it When is an embedding $f:N\to\R^6$ of a 3-manifold $N$ 
isotopic to an embedding $f':N\to\R^6$ such that $f'(N)\subset\R^m$}? 
Here $m=4$ or $m=5$. 
This is a particular case of a classical {\it compression problem} in the 
topology of manifolds 
%studied by Haefliger, Hirsch, Gillman, Tindell, Vrabec, Rourke-Sanderson, 
%Cencelj-Repovs and others in 
[Ha66, Hi66, Gi67, Ti69, Vr89, RS01, Ta04, CR05, KS05, Ta06], \S5. 
This particular case (suggested to the author by Fomenko) is interesting 
because some 3-manifolds appearing in the theory of integrable systems are 
given by a system of algebraic equations implying that the 3-manifold lies 
in $S^2\times S^2\subset\R^5\subset\R^6$ [BF04, Chapter 14]. 

%Embeddings of closed orientable 3-manifolds into $R^6$ are classified by 
%the author in http://arxiv.org/math.GT/0603429 in terms of {\it the Whitney 
%invariant} $W(f)\in H_1(N)$ and {\it the Kreck invariant} 
%$\eta_{W(f)}(f)\in Z/d(W(f))$, where $d(W(f))$ is the divisibility of the 
%projection of $W(f)$ to $H_1(N)/torsion$. 
%The values of the Whitney invariants of embeddings $f:N\to R^6$ such that 
%$f(N)\subset R^5$ are essentially described as follows. 

\smallskip
{\bf Codimension Two Compression Theorem.} 
{\it  Let $N$ be a 3-manifold and $f,f':N\to\R^6$ two embeddings such that 
$f(N)\cup f'(N)$ is contained in either $S^4$ or $S^2\times S^2$ somehow 
embedded into $\R^6$. 
If $H_1(N)$ has no 2-torsion or $W(f)=W(f')$, then $f$ is isotopic to $f'$.}

\smallskip
The proof is given in \S3. 
The proof works for some 4-manifolds different from $S^4$ and $S^2\times S^2$ 
(but not e.g. for $\C P^2\#\overline{\C P^2}$). 

The Codimension Two Compression Theorem implies that when $H_1(N)$ has no 
2-torsion, given an embedding $f':N\to\R^6$ such that $f'(N)\subset S^4$, an 
embedding $f:N\to\R^6$ is compressible to $S^4$ (i.e. is isotopic to an 
embedding with the image in $S^4$) if and only if $f$ is isotopic to $f'$. 
(The same holds for $S^4$ replaced by $S^2\times S^2$).
It would be interesting to obtain a criterion for compressibility into $S^4$ 
not involving another embedding $f'$. 

\smallskip
An embedding $f:N\to\R^6$ is called {\it compressible}, if it is isotopic to an 
embedding $f':N\to\R^6$ such that $f'(N)\subset\R^5$. 
%whose image lies in $\R^5$. 

The Compression Theorem which follows describes compressible embeddings $N\to\R^6$ (or, in the 
equivalent formulation of \S3, the values of the Whitney and the Kreck invariants of 
compressible embeddings). 

We need some definitions. 
Recall that {\it the Rokhlin invariant} $\mu(s)$ of a spin structure $s$ on a 
3-manifold $N$ is $\sigma(V)/8\mod2$, where $V$ is a spin 4-manifold whose spin 
boundary is $(N,s)$. (This is well-defined by the Rokhlin signature theorem.)

An element $u\in H_1(N)$ of order $\le2$ is called {\it spin simple} if 
$\mu(s)=\mu(s')$ for each pair of spin structures $s,s'$ whose difference is 
in $\beta^{-1}u$. 
Here $\beta:H_2(N;\Z_2)\to H_1(N)$ is the (homology) Bockstein homomorphism. 
%$\rho$ is the reduction modulo 2. 

Each order $\le2$ class in 1-dimensional homology is spin simple for:  

$\bullet$
a rational homology 3-sphere (i.e. a 3-manifold 
$N$ such that $H_2(N)=0$), because $\beta$ is injective,  

$\bullet$
$S^1\times S^2$ (because the two different spin structures 
on $S^1\times S^2$ bound spin 4-manifolds homeomorphic to $D^2\times S^2$ and 
$S^1\times D^3$, whose signatures are 0),  

$\bullet$
a connected sum of 3-manifolds with this property. 

An order $\le2$ class in $H_1(S^1\times S^1\times S^1)$ is 0. 
This class is not spin simple [Ki89, V]. 
More generally, for a 3-manifold $N$ without 2-torsion in homology and with 
non-zero intersection form $H_2(N;\Z_2)^3\to\Z_2$ the class $0\in H_1(N)$ 
(i.e. the only order $\le2$ class) is not spin simple [Ka79, Theorem 6.12]. 
We conjecture that for a 3-manifold $N$ with an {\it odd} intersection form 
$H_2(N)^3\to\Z$ (i.e. there exist $x,y,z\in H_2(N)$ such that 
$x\cap y\cap z$ is odd) there is an order $\le2$ class in $H_1(N)$ which is not spin 
simple. 

%THEOREM 6.12. Suppose $N$ is a closed, connected, oriented 3--manifold
%which bounds a framed manifold of signature $k$.  Suppose further that
%there exist elements $x, y, z \in H^1 (N)$ such that $x \cup y \cup z$
%is an odd multiple of the generator of $H^3 (N)$.  Then $N$ bounds a
%framed manifold of signature $k + 8$

\smallskip
{\bf Compression Theorem.} 
{\it An embedding $f:N\to\R^6$ is compressible if and only if 
$2W(f)=0$ and 
either $W(f)$ is not spin simple or $\eta(f,f_0)$ is even for some compressible 
embedding $f_0:N\to\R^6$ such that $W(f_0)=W(f)$ 
(such an embedding $f_0$ exists by an equivalent formulation in \S3).}  

\smallskip
(If $2W(f)=0$, then $d(W(f))=0$ and so $\eta(f,f_0)$ is an integer.)  

The new part of the Compression Theorem is the 'if' part and the necessity of 
the second condition of the 'only if' part. 
(The case $N=S^3$ or $N$ an integral homology sphere of the Compression Theorem 
is known [Ha66, Ta06]; we reprove it for completeness. 
The necessity of $2W(f)=0$ in the `only if' part is easy, see [Vr89] or 
beginning of the second subsection of \S3.)  
See an equivalent formulation of the Compression Theorem in \S3, where we also 
explain more precisely which parts are new. 

The proof of the Compression Theorem (\S3) is based on the Kreck 
Invariant Lemma (this proof does not work for $(4k-1)$-manifolds in $\R^{6k}$). 

\smallskip
{\bf Corollary.} 
{\it (a) If $f:N\to\R^6$ and $l:S^3\to\R^6$ are embeddings, $f$ is 
compressible and $N$ has no 2-torsion in homology and non-zero intersection 
form $H_2(N;\Z_2)^3\to\Z_2$ (the intersection form condition on $N$ can be 
weakened to `$W(f)$ is not spin simple'), then $f\#l$ is compressible (although 
$l$ could be non-compressible).

(b) If two embeddings $N\to\R^5$ are regular homotopic, then their 
compositions with the inclusion $\R^5\to\R^6$ are isotopic.}

\smallskip
Corollary (a) follows from the Compression Theorem, the Addendum to the 
Classification Theorem and [Ka79, Theorem 6.12]. 

The converse to Corollary (b) is trivial by Hirsch-Smale theory. 
Corollary (b) follows by the Classification Theorem and the Wu 
classification of immersions [SST02, Theorems 3.1 and 5.6] because the Whitney 
and the Kreck invariants coincide with regular homotopy invariants $c(f)$ and 
$i_a(f)+\frac32\alpha(N)$ defined in [SST02, Definitions 3.3, 5.1 and 5.3]
(by the proof of the Relation Lemma (b) in \S3, the Compressed Kreck Invariant 
Lemma of \S3 and because each embedding $N\to\R^5$ has a Seifert surface).

%In this talk I shall present explicit examples of embeddings, give 
%definitions of the invariants, state the main result and sketch its proof. 

\bigskip
{\bf Constructions of embeddings. }

%\smallskip
{\it A generalization of Corollary (a) and its proof.}
Cf. [Hu63, Sk07, Hudson Torus Example 3.5, Sk06, Definition of $\bar\mu$ in 
\S5].
For $u\in\Z$ instead of an embedded 3-sphere $2\partial D^4\times x$ we can
take $u$ copies $(1+\frac1n)\partial D^4\times x$ ($n=1,\dots,u$) of 3-sphere 
outside $D^4\times S^2$ `parallel' to $\partial D^4\times x$.
Then we join these spheres by tubes so that the homotopy class of the resulting
embedding $S^3\to S^6-D^4\times S^2\simeq S^6-S^2\simeq S^3$ will be 
$u\in\pi_3(S^3)\cong\Z$.
Let $f$ be the connected sum of this embedding with the standard embedding
$\partial D^2\times S^2\subset\R^6$.

Clearly, $W(f)=u$.
Hence by the Classification Theorem and the Addendum {\it for the generator $t$ of
$\Emb^6(S^3)$ the embedding $f\#kt$ is isotopic to $f$ if and only if $k$ is
divisible by $u$}.
\qed

\smallskip
{\it Any isotopy class of embeddings $N\to\R^6$ can be constructed from a fixed 
embedding $g:N\to\R^6$ by connected summation with an embedding $S^3\to\R^6$} 
linked {\it with the given embedding.}
This follows from the proof of the PL Classification Theorem [Vr77], 
cf. %[Sk06, 
\S5. 
%but not because the restrictions onto $N_0$ of embeddings $N\to\R^6$
%are isotopic to each other [HH63, 3.1.b] - not so for R^5 
Below we present an explicit construction of such an embedding $S^3\to\R^6$
(which generalizes the construction of $f$ from Corollary (a).) 
Note that to construct such an embedding $S^3\to\R^6$ {\it explicitly} (but not 
using embedded surgery or Whitney trick) is an interesting problem.
See some more constructions and remarks in %[Sk06, 
\S5.

\smallskip
{\it Construction of an arbitrary embedding $N\to\R^6$ from a fixed embedding 
$g:N\to\R^5$.}
Represent given $u\in H_1(N)$ by an embedding $u:S^1\to N$. 
Then $\nu_g^{-1}u(S^1)$ is a 3-cycle in $C_g$. 
Recall that any orientable bundle over $S^1$ is trivial. 
Take a section $\overline u:S^1\to \nu_g^{-1}u(S^1)$ of $\nu_g$ 
(e.g. the section pointing from $\R^5$ to $\R^6_+$). 
Take a vector field on $\overline u(S^1)$ normal to $\nu_g^{-1}u(S^1)$ 
(e.g. the upward-looking vector field orthogonal to $\R^5$). 
Extend $\overline u$ along this vector field to a smooth map 
$\overline u:D^2\to S^6$. 
By general position we may assume that $\overline u$ is an embedding
and $\overline u(\Int D^2)$ misses $g(N)\cup \nu_g^{-1}u(S^1)$. 
Take a framing on $\overline u(S^1)$ in $\nu_g^{-1}u(S^1)$. 
Since $\pi_1(V_{4,2})=0$, this framing extends to a 2-framing on 
$\overline u(D^2)$ in $\R^6$. 
Thus $\overline u$ extends to an embedding $\widehat u:D^2\times D^2\to C_g$ 
such that $\widehat u(\partial D^2\times D^2)\subset \nu_g^{-1}u(S^1)$. 
Let 
$$Z:=\nu_g^{-1}u(S^1)-\widehat u(\partial D^2\times\Int D^2)
\bigcup\limits_{\widehat u(\partial D^2\times\partial D^2)}
\widehat u(D^2\times\partial D^2).$$ 
By the Normal Bundle Lemma (a) of \S2 
$\nu_g^{-1}u(S^1)\cong S^1\times S^2$ fiberwise over $S^1$. 
Therefore $Z\cong S^3$. 
Let $f=g\# Z$. 

We have that $\nu_g^{-1}u(S^1)$ spans $\overline{\nu_g}^{-1}u(S^1)$ and 
$g(N)\cap\overline{\nu_g}^{-1}u(S^1)=g(u(S^1))$, where 
$\overline{\nu_g}:S^6-\Int C_g\to N$ is the normal bundle of $g$.  
Since with appropriate orientations $Z$ and $\nu_g^{-1}u(S^1)$ are homologous 
in $C_g$, it follows that $W(f)=u$. 

Any embedding $f':N\to\R^6$ such that $W(f')=u$ can be obtained from $f$ by 
(unlinked) connected summations with the Haefliger trefoil knot [Ha62, 4.1], 
cf. 
%[Sk06, 
\S5, (4).  

\bigskip
{\bf Acknowledgements.} The Classification Theorem was presented at the Yu. P. 
Solovyev memorial conference (Moscow, 2005) and the Conference `Manifolds and 
Their Mappings' (Siegen, 2005). 
The Compression Theorem was presented at the A. B. Sossinskiy mini-conference 
(Moscow, 2007). 
%Cf. [Sk06]. 
I would like to acknowledge D. Crowley, A. T. Fomenko, M. Kreck, 
A. S. Mischenko, N. Saveliev, M. Skopenkov, M. Takase, Su Yang, A. Zhubr 
and the anonymous referees for useful comments.

%\newpage
%\bigskip
\head 2. A surgery proof of the Classification Theorem \endhead

%Recall that $N$ is a closed connected orientable 3-manifold.
%In \S2 by $f:N\to\R^6$ we denote an arbitrary smooth embedding.

%\bigskip 
{\bf The case $N=S^3$.}

The following argument is not used in the proof of the Classification Theorem 
but it is useful to understand why the Kreck invariant is injective.

For $N=S^3$ the Whitney invariant vanishes and the Kreck invariant 
admits the following equivalent definition. 
Let $f:S^3\to\R^6\subset S^6$ be an embedding.
Take a framing $\varphi$ of $f$.
Take the 6-manifold $M=M_\varphi$ obtained from $S^6$ by surgery on this framed 
embedding.
By Alexander duality there is a generator $y\in H_4(M)\cong\Z$ 
such that the intersection of $y$ with the oriented 2-sphere standardly 
linked with $f$ is $+1$.
Since $H_4(M)\cong[M,\C P^\infty]$, there is a closed connected oriented  
4-submanifold $X\subset M$ such that $[X]=2y$. 
Then $\eta(f,g)=\sigma(X)/16$, where $g:S^3\to\R^6$ is the standard 
embedding. 

Clearly, this defines a map $\eta:\Emb^6(S^3)\to\Z$.
This map is a homomorphism (this is a particular case of the second equality 
from the Addendum to the Classification Theorem). 
The Kreck invariant coincides with the Haefliger invariant by [Wa66] or 
[GM86, Remarks to the four articles of Rokhlin, II.2.7 and III.excercises.IV.3, 
Ta04]. 
This map is surjective by [Ha66] or [Ta04]. 

\smallskip
{\it A new proof of the injectivity of $\eta$.}
Take an embedding $f:S^3\to\R^6$ such that $\eta(f)=\eta(f,g)=0$.
Take an orientation-preserving fiberwise diffeomorphism
$\varphi:\partial C_f\to\partial C_g$.
Then $M=M_\varphi\cong C_f\cup_\varphi(-C_g)$.
The manifold $M$ is spin (by the Spin Lemma or because for $N=S^3$ the 
obstructions to extending a spin structure on $C_f$ to that on $M$ assume 
values in zero groups). 
Clearly, $y$ is a joint homology Seifert surface for $f$ and $g_0$. 
By [Wa66] (or by the Twisting Lemma ($\varphi$) below) we can take a framing 
$\varphi$ of $f$ so that $yPDp_1(M)=0$.
Hence $y^3=6\eta(f)=0$.
Hence by the Reduction Lemma in dimension 6 below $f$ is isotopic to $g$.
Thus $\eta$ is injective.
\qed

\bigskip 
{\bf An outline of the proof of the Classification Theorem.}

\smallskip
{\bf Additivity Lemma.}
{\it If $f,f',f'':N\to\R^6$ are embeddings with the same Whitney invariant, 
then $\eta(f,f')+\eta(f',f'')=\eta(f,f'')$.}

\smallskip
{\it Proof.} 
Let $\varphi':\partial C_f\to\partial C_{f'}$ and 
$\varphi'':\partial C_f\to\partial C_{f''}$ be spin isomorphisms. 
Then $\varphi''(\varphi')^{-1}:\partial C_{f'}\to\partial C_{f''}$
is a spin isomorphism. 

Since $H_4(C_f,\partial C_f)\cong[C_f,\C P^\infty]$, there is a proper 
connected oriented 4-submanifold $X\subset C_f$ representing the class 
$2A_f$. 
Analogously by the Agreement Lemma there are proper connected oriented 
4-submanifolds $X'\subset C_{f'}$ and $X''\subset C_{f''}$ representing 
classes $2A_{f'}$ and $2A_{f''}$ such that $\varphi'\partial X=\partial X'$ and 
$\varphi''\partial X=\partial X''$. 
Then the classes of 
$$X\cup(-X')\subset M_{\varphi'},\quad X\cup(-X'')\subset M_{\varphi''}
\quad\text{and}\quad X'\cup(-X'')\subset M_{\varphi''(\varphi')^{-1}}$$ 
are twice the joint homology Seifert surfaces for $f$ and $f'$, $f$ and $f''$, 
$f'$ and $f''$, respectively.  
(In this proof $\cup$ means the union along common boundary.)
Hence by the Novikov-Rokhlin additivity
$$16\eta(f,f')+16\eta(f',f'')=\sigma(X\cup(-X'))+\sigma(X'\cup(-X''))=
\sigma(X\cup(-X''))=16\eta(f,f'').\ \qed$$

%\smallskip
{\it Proof of the Classification Theorem.}
The surjectivity of $W$ is proved at the end of \S1. 
So it remains to prove the bijectivity of $\eta_u$ for each $u\in H_1(N)$.
By the Addendum $\eta_u$ is surjective.
The injectivity of $\eta_u$ is implied by the Additivity Lemma and the following 
Injectivity Lemma.  
\qed

\smallskip
{\bf Injectivity Lemma.}
{\it If $f,f':N\to\R^6$ are embeddings with the same Whitney invariant and 
$\eta(f,f')\equiv0\mod d(W(f))$, then $f$ is isotopic to $f'$.}

\smallskip
For the proof we need the following two results (which are proved below 
in this section). 

\smallskip
{\bf Reduction Lemma in dimension 6.}
{\it Two embeddings $f,f':N\to\R^6$ are isotopic if and only if there is a spin 
isomorphism $\varphi:\partial C_f\to\partial C_{f'}$ and a joint homology 
Seifert surface $y\in H_4(M_\varphi)$ for $f$ and $f'$ such that 
$yPDp_1(M_\varphi)=y^3=0$.}

\smallskip
{\bf Twisting Lemma.} {\it In the notation from the 
definition of the Kreck invariant one can change 

(y) the joint homology Seifert surface $y$ so that $\sigma_{2y}(M_\varphi)/16$ 
would change by adding $d(W(f))$.

($\varphi$) the spin isomorphism $\varphi:\partial C_f\to \partial C_{f'}$ over 
$B^3$ (i.e. over the top cell of $N$) and the joint homology Seifert surface 
$y$ so that $y^3$ and $yPDp_1(M_\varphi)$ would change by adding 1 and 4, 
respectively.}  

\smallskip
{\it Proof of the Injectivity Lemma.} 
Take a spin isomorphism $\varphi$. 
By the Agreement Lemma there is a joint homology Seifert surface $y$ for $f$ and 
$f'$. 
Then 
$$0\equiv\frac{\sigma_{2y}(M_\varphi)}{16}\equiv
\frac{yPDp_1(M_\varphi)-4y^3}{24}\mod d(W(f)).$$ 
Change $y$ by the Twisting Lemma ($y$) so that $yPDp_1(M_\varphi)-4y^3=0$. 
Change $\varphi$ and $y$ by the Twisting Lemma ($\varphi$) so that 
$yPDp_1(M_\varphi)=4y^3=0$. 
Now the Injectivity Lemma follows from the Spin Lemma and the Reduction Lemma 
in dimension 6. 
\qed

\bigskip 
{\bf Proof of the Agreement Lemma.}

By a section we always mean a section $N\to\partial C_f$ of 
$\nu_f:\partial C_f\to N$.
For a map $\xi:P\to Q$ between a $p$- and  a $q$-manifold denote by 
$$\xi^!:=PD\circ\xi^*\circ PD:H_i(Q,\partial Q)\to H_{p-q+i}(P,\partial P)$$ 
the `preimage' homomorphism.  
Denote by $\overline e(f)\in\Z$ and $\overline w_i(N)\in H^i(N;\Z_2)$ the normal 
Euler and Stiefel-Whitney classes [MS74]. 

\smallskip
{\bf Normal Bundle Lemma.}
{\it (a) The normal bundle of $f$ is trivial. 

(b) Every section on $N_0$ extends to that of $N$.} 

\smallskip
{\it Proof of (a).} 
Since $N$ is a closed connected orientable 3-manifold, we have 
$\overline w_1(N)=0$, $\overline w_2(N)=0$ and $\overline e(f)=0$.
Hence part (a) follows by [DW59]. 
\qed

%, i.e any embedding $f:N\to\R^6$ has a framing. A framing of $f$ defines an 
%isomorphism $\partial C_f\cong N\times S^2$ and the K\"unneth isomorphism 
%$k_\xi:H_3(\partial C_f)\overset\xi\to\to H_3(N)\oplus H_1(N)$.
%The latter depends only on the section $\xi:N_0\to\partial C_f$ formed by the 
%first vectors of the framing. 

\smallskip
{\it Proof of (b).} 
By (a), sections over any subcomplex $X\subset N$ are in 1--1 correspondence 
with $[X;S^2]$. 
Consider the following diagram
$$\CD [N;S^2] @>> \deg > H_1(N) \\ 
@VV r V  @VV \cong V\\
[N_0;S^2] @>> \deg > H_1(N_0,\partial N_0)\endCD.$$ 
Here $r$ is the restriction map, $\cong$ is the inclusion-induced isomorphism 
and $\deg[g]:=g^!(*)$.    
Now part (b) follows because the lower $\deg$ is bijective and the upper $\deg$ 
is surjective. 
\qed

\smallskip
Denote by 

$\bullet$
$d(\xi,\zeta)\in H_1(N)$ the difference element of sections $\xi,\zeta$,  

$\bullet$
$\xi^\perp$ the orthogonal complement to a section $\xi:N\to\partial C_f$ in 
the normal bundle of $f$ ($\xi^\perp$ is an oriented $S^1$-bundle), 

$\bullet$
$A_{f,\xi}$ the image of $[N]$ under the composition 
$$H_3(N)\overset{\xi_*}\to\to H_3(\partial C_f)\overset{in_*}\to\to H_3(C_f)
\overset{AD}\to\to H^2(N)\overset{PD}\to\to H_1(N).$$ 
%of Alexander duality and Poincar\'e isomorphisms. 
The sign in the following formulas (which is denoted by $\pm$) is fixed, i.e. 
does not depend on $N$, $f$, $\xi$ etc. 

\smallskip
{\bf Difference Lemma.} 
{\it (a) For sections $\xi$ and $\zeta$ we have 
$\pm d(\xi,\zeta)=A_{f,\xi}-A_{f,\zeta}$. 
 
(b) For sections $\xi$ and $\zeta$ we have 
$PDe(\xi^\perp)-PDe(\zeta^\perp)=\pm2d(\xi,\zeta)$. 
 
(c) If $f=f'$ on $N_0$ and $\xi,\xi'$ are sections of $f,f'$ such that 
$\xi=\xi'$ on $N_0$, then 
$W(f)-W(f')=A_{f,\xi}-A_{f',\xi'}$.}

\smallskip
Part (a) follows by Alexander duality, cf. 
[BH70, Lemme 1.2, Bo71, Lemme 3.2.a]. 
The proof of part (b) is analogous to [BH70, Lemme 1.7, Bo71, Lemme 3.2.b]. 

\smallskip
{\it Proof of (c).} 
Denote by $A_0$ the composition 
$H_3(S^6-f(N_0))\to H^2(N_0)\to H^2(N)\to H_1(N)$  
of the Alexander duality, the restriction and Poincar\'e isomorphisms.
Then  
$$A_{f,\xi}-A_{f',\xi'}=A_0[\xi|_{B^3}\cup\xi'|_{\overline B^3}]=
A_0[f|_{B^3}\cup f'|_{\overline B^3}]=W(f)-W(f').$$
Here $\overline B^3$ is $B^3$ with reversed orientation, 

$\bullet$
the first equality follows because $\xi=\xi'$ over $N_0$, 

$\bullet$
the second equality follows because in $S^6-f(N_0)$ we can shift $\xi|_{B^3}$ 
and $\xi'|_{B^3}$ to $f|_{B^3}$ and $f'|_{B^3}$, 

$\bullet$
the third equality follows by Alexander duality. 
\qed

\smallskip
A section $\xi:N\to\partial C_f$ is called {\it unlinked} if 
$A_{f,\xi}=0$ [BH70].
  
%$i_*\xi_*[N]=0\in H_3(C_f)$, 

\smallskip
{\bf Unlinked Section Lemma.}
{\it (a) An unlinked section exists and is unique on $N_0$ up to fiberwise 
homotopy, cf. [HH63, 4.3, BH70, Proposition 1.3 and Lemme 1.2].

(b) If $W(f)=W(f')$, then any spin isomorphism maps an unlinked section of $f$ 
to an unlinked section of $f'$, cf. [KS05]. 

(c) If $\xi$ is an unlinked section, then $\xi_*[N]=\partial A_f$.}  

\smallskip
{\it Proof of (a).} For fixed $\xi$ the map $\zeta\mapsto d(\xi,\zeta)$ 
defines a 1--1 correspondence between sections over $N_0$ and 
$H_1(N_0,\partial N_0)\cong H_1(N)$. 
This and the Difference Lemma (a) imply the uniqueness. 
This, the Difference Lemma (a) and the Normal Bundle Lemma (b) imply the 
existence. 
\qed

\smallskip
{\it Proof of (b).} Since the assertion is invariant under isotopy of $f$, we 
may assume that $f=f'$ on $N_0$ and $\varphi$ is the identity on $N_0$.
Let $\xi$ be the unlinked section of $f$. 
Apply the Difference Lemma (c) to $\xi'=\varphi\xi$. 
We obtain $A_{f',\varphi_*\xi_*}=0$, i.e. the section $\varphi\xi$ of $f'$ 
is unlinked.  
\qed

\smallskip
{\it Proof of (c).} 
Since $\xi$ is unlinked and $H_4(C_f,\partial C_f)$ is generated by $A_f$, 
it follows that $\xi_*[N]=k\partial A_f$ for some integer $k$. 
We have $k=1$ because $[N]=\nu_*\xi_*[N]=k\nu_*\partial A_f=k[N]$. 
Here the first equality holds because $\xi$ is a section and the last by the 
following lemma. 
\qed

\smallskip
{\bf Alexander Duality Lemma.} {\it The composition $\nu_*\partial$ equals to 
the composition of the Alexander duality and Poincar\'e isomorphisms.} 

\smallskip
{\it Proof.} The composition $\nu_*\partial$ goes under the inclusion 
$J:(C_f,\partial C_f)\to(S^6,S^6-\Int C_f)$ to the composition 
$$H_4(S^6,S^6-\Int C_f)\overset{\overline\partial}\to\to H_3(S^6-\Int C_f)
\overset{\overline\nu_*}\to\to H_3(N)$$ (of the map $\overline\nu_*$ 
induced by 
the normal vector bundle $\overline\nu:S^6-\Int C_f\to N$ and the boundary 
isomorphism $\overline\partial$). 
Now the Lemma follows because 
$\overline\nu_*\overline\partial J_*=\nu_*\partial$. 
\qed

\smallskip
{\it Proof of the Agreement Lemma.} 
Consider the following fragment of the Gysin sequence for the (trivial) bundle 
$\nu:=\nu_f$: 
$$0\to H_1(N)\overset{\nu^!}\to\to H_3(\partial C_f)\overset{\nu_*}\to
\to H_3(N)\to0.$$
We see that for each section $\xi:N\to\partial C_f$ the map 
$$\nu_*\oplus\xi^!:H_3(\partial C_f)\to H_3(N)\oplus H_1(N)$$ 
is an isomorphism. 
By the Alexander Duality Lemma 
$\nu_{f,*}\partial A_f=[N]=\nu_{f',*}\partial A_{f'}$. 
By the Unlinked Section Lemma (a) there exists an unlinked section $\xi$. 
By the Unlinked Section Lemma (c) and since the normal bundle of 
$\xi:N\to\partial C_f$ is isomorphic to $\xi^\perp$, we have 
$\xi^!\partial A_f=\xi^!\xi_*[N]=PDe(\xi^\perp)$. 
Hence by the Unlinked Section Lemma (b) 
$(\varphi\xi)^!\partial A_f=PDe((\varphi\xi)^\perp)=PDe(\xi^\perp)=
\xi^!\partial A_f$. 

Therefore $\varphi_*\partial A_f=\partial A_{f'}$. 
Now the Lemma follows by looking at the segment of (Poincar\'e dual of) the 
Mayer-Vietoris sequence: 
$$H_4(\partial C_f)\to H_4(M_\varphi)
\to H_4(C_f,\partial C_f)\oplus H_4(C_{f'},\partial C_{f'})
\to H_3(\partial C_f)\quad\qed$$

\smallskip
%\bigskip
{\bf Proof of the Independence and Twisting Lemmas.}

Note that in the Independence Lemma $y$ can be changed without changing 
$\varphi$, but a change of $\varphi$ forces a change of $M=M_\varphi$ and hence 
one of $y$. 

\smallskip
{\it Proof that $\sigma_{2y}(M)/16\mod d(W(f))$ is independent of the choice of 
joint homology Seifert surface $y$ when $\varphi$ is fixed.}
From the Mayer-Vietoris sequence at the end of the proof of the Agreement Lemma 
we see that the choice of $y$ is in adding a class $y_1\in H_4(M)$ coming from 
$H_4(\partial C_f)$. 

By the Normal Bundle Lemma (a) there exists a framing of $\nu=\nu_f$. 
This framing defines a (fiberwise over $N$) diffeomorphism 
$\partial C_f\cong S^2\times N$. 
Identify $\partial C_f$ with  $S^2\times N$ by this diffeomorphism. 
This diffeomorphism induces isomorphisms 
$$H_4(\partial C_f)\cong H_2(N)\quad\text{and}\quad
\nu_*\oplus\xi^!:H_3(\partial C_f)\to H_3(N)\oplus H_1(N),$$  
where $\xi$ is the section formed by the first vectors of the framing. 

By the first isomorphism we may assume that $y_1=[S^2\times L]$ for a certain 
sphere with handles (i.e. closed orientable surface) $L\subset N$. 
We have 
$$(y+y_1)PDp_1(M)-yPDp_1(M)=y_1PDp_1(M)=p_1(S^2\times L)=3\sigma(S^2\times L)=
0.$$
Here the second equality follows because the normal bundles of $L\subset N$ 
and $\partial C_f\subset C_f$ are trivial, hence that of $S^2\times L\subset M$ 
is trivial. 

Let $Z\subset N$ be an oriented circle representing $W(f)$ and $*\in S^2$ 
a point. 
Denote by $\cap$ the intersection in $\partial C_f$. 
Then
$$\frac{4(y+y_1)^3-4y^3}{24}\overset{(1)}\to=\frac{y^2y_1}2\overset{(2)}\to=
\frac{\partial A_f\cap\partial A_f\cap y_1}2\overset{(3)}\to=
\pm[*\times Z]\cap[S^2\times L]\overset{(4)}\to=\pm W(f)\bigcap\limits_N[L]$$ 
is divisible by $d(W(f))$. 
So $\sigma_{2y+2y_1}(M)-\sigma_{2y}(M)$ is divisible by $16d(W(f))$. 
 
Here (1) holds because the normal bundle of 
$\partial C_f\subset C_f$ is trivial, so $y_1^2=0$, 

(2) holds because $\partial A_f=\partial(y\cap C_f)=y\cap\partial C_f$,   

(3) holds because by the Alexander Duality Lemma $\nu_*\partial A_f=[N]$, by 
the Framing Lemma below $2\xi^!\partial A_f=\pm2W(f)$, hence 
$$2\partial A_f=2[*\times N]\pm2[S^2\times Z],\quad\text{so}\quad
2\partial A_f\cap2\partial A_f=\pm8[*\times Z].$$ 

(4) holds because $[Z]=W(f)$. 
\qed

\smallskip
{\bf Framing Lemma.} {\it If a section $\xi$ of $\nu_f$ extends to a framing of 
$\nu_f$, then $2\xi^!\partial A_f=\pm2W(f)$.}

\smallskip
{\it Proof.} Let $\zeta$ be an unlinked section for $f$. 
Then up to sign we have 
$$2\xi^!\partial A_f\overset{(1)}\to=2\xi^!\zeta_*[N]\overset{(2)}\to=
2d(\xi,\zeta)\overset{(3)}\to=PDe(\xi^\perp)-PDe(\zeta^\perp)\overset{(4)}\to=
PDe(\zeta^\perp)\overset{(5)}\to=2W(f).$$ 
Here  (1) holds by the Unlinked Section Lemma (c), 

(2) and (4) hold because $\xi$ extends to a framing, 

(3) holds by the Difference Lemma (b), 

(5) holds by the Boechat-Haefliger Invariant Lemma below. 
\qed

\smallskip
Take an unlinked section $\xi:N\to\partial C_f$ and define 
{\it the Bo\'echat-Haefliger invariant}
$$BH:\Emb\phantom{}^6(N)\to H_1(N)\quad\text{by}\quad BH(f):=PDe(\xi^\perp).$$
We have $BH(f)\mod2=PDw_2(N)=0$.
Note that this invariant unlike the Whitney invariant is independent of the
choice of the embedding $g$ (but just as the Whitney invariant depends on the 
choice of an orientation of $N$).

\smallskip
{\bf The Bo\'echat-Haefliger Invariant Lemma.} $BH(f)=\pm2W(f)$.

\smallskip
{\it Proof.} The section normal to $\R^5\subset\R^6$ is unlinked for the 
embedding $g$ fixed in the definition of the Whitney invariant. 
Since every codimension 2 embedding into Euclidean space of an orientable 
manifold has a trivial normal bundle, we have $BH(g)=0$. 
Thus it suffices to prove that $BH(f)-BH(f')=\pm2(W(f)-W(f'))$. 

Since both $BH(f)$ and $W(f)$ are invariant under isotopy of $f$, we may assume 
that $f=f'$ on $N_0$.
Let $\xi':N\to\partial C_{f'}$ be an unlinked section for $f'$. 
Then $\xi'|_{N_0}$ is a section for $f$ on $N_0$.
Extend $\xi'|_{N_0}$ to a section $\xi:N\to C_f$ of $f$.
Take an unlinked section $\zeta$ for $f$. 
Then up to sign we have 
$$BH(f)-BH(f')\overset{(1)}\to=
PDe(\zeta^\perp)-PDe(\xi^\perp)\overset{(2)}\to=
2d(\zeta,\xi)\overset{(3)}\to=
2A_{f,\zeta}-A_{f,\xi}\overset{(4)}\to=$$ 
$$=2A_{f,\xi}\overset{(5)}\to=
2A_{f,\xi}-2A_{f',\xi'}\overset{(6)}\to=
2W(f)-2W(f').$$
Here (1) holds because $\xi=\xi'$ on $N_0$; 

(2) holds by the Difference Lemma (b); 

(3) holds by the Difference Lemma (a); 

(4) holds because $\zeta$ is unlinked; 

(5) holds because $\xi'$ is unlinked; 

(6) holds by the Difference Lemma (c). 
\qed

\smallskip
Note that $\xi^!\partial A_f=\pm2W(f)$ for an unlinked section $\xi$ by the 
proof of the Agreement Lemma and the Bo\'echat-Haefliger Invariant Lemma, cf. 
[KS05, Euler Class Lemma]. 

\smallskip
{\it Proof of the Twisting Lemma (y).} Analogously to the above proof of the 
independence on $y$ 
because Poincar\'e duality implies that by changing $y_1$ and thus $L$ 
we may obtain as $W(f)\bigcap\limits_N[L]$ any integer divisible by $d(W(f))$. 
\qed

\smallskip
{\bf Spin Framing Lemma.}
{\it A spin framing $\varphi$ is unique over $N_0$ up to fiberwise isotopy.} 

\smallskip
{\it Proof.} 
It suffices to prove that if framed embeddings $\xi,\xi_1:N_0\times D^2\to\R^6$ 
are isotopic, then they are isotopic relative to $N_0\times0$.

By general position for large enough $m$ the compositions of $\xi,\xi_1$ with 
the inclusion $\R^6\subset\R^m$ are isotopic relative to $N_0$.
Hence the difference $d(\xi,\xi_1)\in[N_0,SO]$ is zero.
The map $[N_0,SO_3]\to[N_0,SO]$ induced by the 
inclusion $SO_3\subset SO$ is a 1--1 correspondence (indeed, $N_0$ retracts to 
a 2-dimensional complex, the group $\pi_1(SO_3)$ is stable and the group 
$\pi_2(SO_3)=\pi_2(SO_4)=0$ is stable).
Thus $\xi$ and $\xi_1$ are isotopic relative to $N_0\times0$.
\qed

\smallskip
Note that the above proof together with the classification of spin structures 
on $N$ imply that $\varphi|_{N_0}$ is uniquely defined (up to fiberwise 
homotopy) by the condition that $M=M_\varphi$ is spin. 

\smallskip
{\it Proof of the Independence Lemma.} 
Since $\sigma_{2y}(M)/16\mod d(W(f))$ is independent of a change of $y$, it is 
independent of a change of $\varphi$ by a fiberwise isotopy (for a certain 
corresponding change of $y$).
By the Spin Framing Lemma this residue is independent of a change of 
$\varphi|_{N_0}$ (for a certain corresponding change of $y$).

Since the normal bundle of $f$ is trivial, we have a homeomorphism 
$S^2\times N\to \partial C_f$ depending on a framing. 
However, the images of $S^2\times B^3$ and $S^2\times N_0$ do not depend on the 
framing. 
So we identify these images with $S^2\times B^3$ and $S^2\times N_0$, 
respectively. 
Make the same identification for $f'$. 
Now assume that $\varphi=\varphi'$ over $S^2\times N_0$ and that 
$\varphi|_{S^2\times B^3}$ is obtained from $\varphi'|_{S^2\times B^3}$ by 
twisting with $\alpha\in\pi_3(SO_3)\cong\Z$. 
It suffices to consider the case when $\alpha$ is the generator. 
Denote $a:=[\C P^2]\in H_4(\C P^3)$. 
The Lemma holds because by the following Cobordism Lemma we have 
$$\sigma_{2y}(M_\varphi)-\sigma_{2y'}(M_{\varphi'})=
\sigma_{2a}(\C P^3)=(4a^2\cap2a-(2a)^3)/3=0.\quad\qed$$

{\bf Cobordism Lemma.} {\it In the above notation the pairs $(M_{\varphi'},y')$ 
and $(M_\varphi\sqcup\C P^3,y+[\C P^2])$ are cobordant for certain joint 
homology Seifert surfaces $y'$ and $y$.}

\smallskip
{\it Proof.} 
Define 
$$W:=C_f\times[0,3]\bigcup\limits_{\overline\alpha}\ C_{f'}\times[0,3],\quad
\text{where}$$ 
$$\overline\alpha\ :\ \partial C_f\times[0,3]-S^2\times\Int B^3\times(1,2)\ \to
\ \partial C_{f'}\times[0,3]-S^2\times\Int B^3\times(1,2)$$
$$\text{is defined by}\quad \overline\alpha(x,t):=\cases 
(\alpha(b)s,b,t)&x=(s,b)\in S^2\times B^3\text{ and }t\in[2,3]\\  
(\varphi(x),t)&\text{otherwise}\endcases.$$
By the proof of the Agreement Lemma $\varphi_*\partial A_f=\partial A_{f'}$ 
and the same for $\varphi$ replaced by $\varphi'$. 
Looking at (the Poincar\'e dual of) the Mayer-Vietoris sequence for $W$ (cf. 
the end of the proof of the Agreement Lemma) we see that there is a class 
$$w\in H_5(W)\quad\text{such that}
\quad w|_{C_f\times[0,3]}=p^!A_f\quad\text{and}\quad 
w|_{C_{f'}\times[0,3]}=(p')^!A_{f'},$$ 
where $p:C_f\times[0,3]\to C_f$ and $p':C_{f'}\times[0,3]\to C_{f'}$ are the 
projections. 
We have $\partial W\cong M_\varphi\sqcup(-M_{\varphi'})\sqcup E$, where 
$E$ is the remaining boundary component. 
Clearly, 
$(w|_{M_\varphi})|_{C_f}=A_f$ and the same when either $f$ is replaced 
by $f'$ or $\varphi$ is replaced by $\varphi'$ or both.  
Hence $(W,w)$ is a cobordism between $(M_{\varphi'},y')$ and 
$(M_\varphi\sqcup E,y\sqcup w|_E)$ for certain joint homology Seifert 
surfaces $y'$ and $y$. 
We have 
$$E\cong S^2\times B^3_+\times[1,2]\bigcup\limits_{\widehat\alpha} 
S^2\times B^3_-\times[1,2]\cong\frac{S^2\times B^4}
{\{(s,b)\sim(\alpha(b)s,\sigma b)\}_{(s,b)\in S^2\times B^3_+}}\cong\C P^3,
\quad\text{where}$$
$$\widehat\alpha:
S^2\times\partial(B^3_+\times[1,2])\to S^2\times\partial(B^3_-\times[1,2])
\quad\text{maps $(s,b,t)$ to}\quad
\cases(s,b,t) &t<2 \\ (\alpha(b)s,b,2)& t=2\endcases,$$ 
$\partial B^4_+=B^3_+\cup B^3_-$ and $\sigma:B^3_+\to B^3_-$ is the symmetry 
with respect to $\partial B^3_+=\partial B^3_-$. 
The latter diffeomorphism to $\C P^3$ is well-known and is proved using a 
retraction to the dual $(\C P^1)^*\subset\C P^3$ of the complement to a 
(trivial) tubular neighborhood of $\C P^1\subset\C P^3$. 

From (the Poincar\'e dual of) the Mayer-Vietoris sequence for $E$ we see that 
the intersections of $w$ with the parts of $E$ are represented by 
$*\times B^3_\pm\times[1,2]$. 
(On the intersection of the parts these manifolds do not coinside, but 
the represented classes are homologous.) 
Hence under the second homeomorphism $w|_E$ goes to a class whose intersection 
with $[S^2\times 0]$ is $+1$.
Therefore under the compositions of the diffeomorphisms $w|_E$ goes to 
a class whose intersection with $[\C P^1]$ is $+1$, i.e. to $[\C P^2]$. 
\qed

\smallskip
The Twisting Lemma ($\varphi$) is proved analogously to the above proof. 

%So $p_1(E_\psi)\cup y|_{E_\psi}=4y|_{E_\psi}^3=4\psi$ by [Wa66, Theorem 4]. 
%$E_0=S^2\times S^4$.)
%$\chi_{S^3,1}(S^6)$ is obtained from $S^6$ by surgery on the trivial knot with 
%the framing $q:S^3\times D^3\to S^6$ such that the linking coefficient of 
%$q(S^3\times0)$ and $q(S^3\times z)$ is 1 for $z\ne0$. 
%The second diffeomorphism holds because $E$ is obtained from $S^2\times S^4$ 
%by first deleting $S^2\times B^4_-$ thus obtaining $S^2\times B^4_+$, then 
%representing $\partial B^4_+=B^3_+\cup B^3_-$ and at last identifying 
%$(s,b)\in S^2\times B^3_+$ with $(\alpha(b)s,b)\in S^2\times B^3_-$. 

%\newpage
\bigskip
{\bf The Reduction Lemma.}

A map between connected spaces is called {\it $m$-connected} if it induces an 
isomorphism on $\pi_i$ for $i<m$ and an epimorphism on $\pi_m$.

\smallskip
{\bf Diffeomorphism Theorem.} 
{\it Let $q\ge3$ be odd and $W$ be a compact simply-connected $(2q+1)$-manifold 
(embedded into $\R^{5q}$) whose boundary is the union along the common boundary 
of compact simply-connected $2q$-manifolds $M_0$ and $M_1$ with the same Euler 
characteristic.
Let $\pi:B\to BO$ be a fibration and $\bar\mu:W\to B$ be a lifting of the Gauss
map $\mu:W\to BO$. 

If $\pi_1(B)=0$ and $\bar\mu|_{M_0}$ and $\bar\mu|_{M_1}$ are 
$q$-connected, then the identification $\partial M_0=\partial M_1$ extends 
to a diffeomorphism $M_0\cong M_1$.}

\smallskip
This result follows by [Kr99, Theorems 4 and 5.i] (cf. [Kr99, Corollary 4] 
where the obtained diffeomorphism extends the identification 
$f:\partial M_0\to\partial M_1$). 

A. Zhubr kindly informed me that for $q=3$ Diffeomorphism Theorem can possibly 
be obtained using the version of surgery described in [Zh75, 4.1], and that the 
Decomposition Theorem used in [Zh75, 4.1] has a simpler proof. 

%\smallskip
Denote by $BO\left<m\right>$ the (unique up to homotopy equivalence) 
$(m-1)$-connected space for which there exists a fibration 
$p:BO\left<m\right>\to BO$ inducing an isomorphism on $\pi_i$ for $i\ge m$.
(M. Kreck kindly informed me that in [Kr99, Definition of $k$-connected cover 
on p. 712] $X\left<k\right>$ should be replaced by $X\left<k+1\right>$.)

For a fibration $\pi:B\to BO$ denote by $\Omega_q(B)$ the group of bordism
classes of liftings $\bar\mu:Q\to B$ of the stable Gauss map (i.e. a classifying 
map of the stable normal bundle) $\mu:Q\to BO$, where $Q$ is a $q$-manifold 
embedded into $\R^{3q}$.
This should be denoted by $\Omega_q(\pi)$ but no confusion would arise.
(This group is the same as $\Omega_q(B,\pi^*t)$ in the notation of [Ko88].)

\smallskip
{\bf Reduction Lemma.}
{\it Let $q\ge3$ be odd, $3\le n\le2q-3$ and $N$ a closed connected 
$n$-manifold.
Two embeddings $f,f':N\to\R^{2q}$ are isotopic if for some isomorphism 
$\varphi:\partial C_f\to\partial C_{f'}$ and some embedding 
$M_\varphi\to S^{5q}$ there exist

a space $C$,

a map $h:M_\varphi\to C$ whose restrictions to $C_f$ and to $C_{f'}$ are 
$q$-connected, and

a lifting $l:M_\varphi\to BO\left<q+1\right>$ of the Gauss map 
$\mu:M_\varphi\to BO$ such that 
$$[h\times l]=0\in\Omega_{2q}(C\times BO\left<q+1\right>).$$}
This situation is explained in the following diagram: 
$$\minCDarrowwidth{5pt}\CD
B:=C\times BO\left<q+1\right> @>> \pr_2 > BO\left<q+1\right> \\
@AA h\times l A                   @VV p V          \\
M_\varphi=C_f\cup_\varphi(-C_{f'})                 @>> \mu   > BO
\endCD.$$

%The `only if' direction is obvious. Let us prove the `if' direction.

{\it Proof.} 
The Euler characteristics of $C_f$ and $C_{f'}$ are the same by Alexander
duality.
Denote $B:=C\times BO\left<q+1\right>$.
Since $[h\times l]=0$, it follows that there is a simply-connected zero bordism 
$\bar\mu:W\to B$ of $h\times l$.

Since is $C_f$  simply-connected and $h|_{C_f}$ is $q$-connected, $B$ is 
simply-connected.  
Since $h|_{C_f}$ and $h|_{C_{f'}}$ are $q$-connected and 
$BO\left<q+1\right>$ is $q$-connected, it follows that 
$\bar\mu|_{C_f}=(h\times l)|_{C_f}$ and 
$\bar\mu|_{C_{f'}}=(h\times l)|_{C_{f'}}$ are $q$-connected.
Therefore by the Diffeomorphism Theorem 
$\varphi:\partial C_f\to\partial C_{f'}$ 
extends to an orientation-preserving diffeomorphism $C_f\to C_{f'}$.
Since any orientation-preserving diffeomorphism of $\R^{2q}$ is isotopic to the
identity, it follows that $f$ and $f'$ are isotopic.
\qed

\smallskip
{\it Proof of the Reduction Lemma in dimension 6.} 
Take any embedding $M_\varphi\to S^{16}$. 
Let $C:=\C P^\infty$. 
Let $h_f$ be any map corresponding to $A_f$ under the bijection  
$H_4(C_f,\partial C_f)\to[C_f,\C P^\infty]$. 
Define $h'$ analogously. 
Since $y$ is a joint homology Seifert surface for $f$ and $f'$, it follows that 
$\varphi_*\partial A_f=\partial A_{f'}$, i.e.  
$h_f|_{\partial C_f}$ is homotopic to $h_{f'}|_{\partial C_{f'}}\circ\varphi$. 
Therefore by the Borsuk homotopy extension theorem $h_{f'}$ is homotopic to a
map $h':C_{f'}\to\C P^\infty$ such that $h'\varphi=h_f$ on $\partial C_f$.
Set $h:=h_f\cup h'$. 
Since $BO\left<4\right>=BSpin$, by the Spin Lemma there is a lifting 
$l:M_\varphi\to BSpin$. 

We have  
$\Omega_6(\C P^\infty\times BO\left<4\right>)\cong\Omega^{spin}_6(\C P^\infty)$ 
is the group of spin cobordism classes of maps $y:Q\to\C P^\infty$ from spin 
6-manifolds $Q$.
Define a map
$$\omega:\Omega^{spin}_6(\C P^\infty)\to\Z\oplus\Z\quad\text{by}\quad
\omega(y:Q\to\C P^\infty):=(\widehat yPDp_1(Q),\widehat y^3),$$
where $\widehat y:=y^![\C P^2]$. 
Then {\it $\omega$ is a monomorphism} [Fu94, proof of Proposition 1.1].
Now the Reduction Lemma in dimension 6 follows by the Reduction Lemma for $q=3$. 
\qed

%\newpage
%\bigskip
\head 3. Proof of the Compression Theorems \endhead

{\bf Seifert surfaces and the Kreck invariant.}

A {\it Seifert surface} of an embedding $f:N\to\R^m$ is an extension of $f$ to 
an embedding $\t f:V\to\R^m$ of a connected oriented 4-manifold $V$ with 
boundary $\partial V=N$ (the normal bundle of $\t f$ is not necessarily 
framed).

Each embedding $f:N\to\R^5$ has a Seifert surface $\t f:V\to\R^5$ 
[Ki89, VIII, Theorem 3]. 
%Note that {\it any smooth embedding of a closed orientable $n$-manifold in 
%$\R^{n+2}$ has a Seifert surface} [Ki89, VIII, Theorem 3]. 
%[Br93, p. 515]. 
The following result (which is not used in this paper) was 
communicated to me by P. Akhmetiev as well-known. 
Cf. [Ta04, 3.3, Ta06, Proposition 2.5].  

\smallskip
{\bf Existence Lemma.} 
{\it Each embedding $f:N\to\R^6$ has a Seifert surface.}  

\smallskip
{\it Proof (folklore).} 
By the Unlinked Section Lemma (a) of \S2 there is an unlinked section 
$\xi:N\to\partial C_f$.  
Consider the following diagram: 
$$\CD [C_f;\C P^\infty] @ >>\psi> H_4(C_f,\partial C_f)\supset A_f\\
@VV r V @VV\partial V \\
[\partial C_f;\C P^\infty] @>>\psi> H_3(\partial C_f)\supset \xi_*[N]\endCD.$$ 
Here $r$ is induced by restriction and the maps $\psi$ are bijections. 
In this diagram $\C P^\infty$ can be replaced by $\C P^3$. 
Take a map $\bar h_f:\partial C_f\to\C P^3$ transverse to $\C P^2$ 
and such that $\bar h_f^{-1}(\C P^2)=\xi(N)$ and $\psi[\bar h_f]=\xi_*[N]$.  
By the Unlinked Section Lemma (c) of \S2 there exists an extension 
$h_f:C_f\to\C P^3$ of $\bar h_f$ transverse to $\C P^2$ and such that 
$\psi[h_f]=A_f$. 
Then $h_f^{-1}(\C P^2)$ is an orientable 4-manifold with boundary 
$\xi(N)\subset\partial C_f$. 
The union of $h_f^{-1}(\C P^2)$ and arcs $x\xi(x)$, $x\in N$, can be smoothed 
to give a Seifert surface for $f$. 
\qed

\smallskip
{\bf Seifert Surface Lemma.} 
{\it Let $f,f':N\to\R^6$ be two embeddings such that $W(f)=W(f')$ and let 
$\t f:V\to\R^6$, $\t f':V'\to\R^6$ be Seifert surfaces for $f$ and $f'$.
Denote $V_0:=\t f(V)\cap C_f$ and $V_0':=\t f'(V')\cap C_{f'}$.

(a) There is a spin isomorphism $\varphi:\partial C_f\to\partial C_{f'}$ such 
that $\varphi(\partial V_0)=\partial V_0'$. 
 
(b) For such a spin isomorphism the homology class of 
$Y:=V_0\cup_\varphi(-V_0')$ 
in $H_4(M_\varphi)$ is a joint homology Seifert surface for $f$ and $f'$. 

(c) If both $\t f$ and $\t f'$ have normal sections, then $\overline e$ is in 
the image of the inclusion homomorphism $H_2(\partial V_0)\to H_2(Y)$ and so 
$\overline e\cap\overline e=0$.} 

\smallskip
{\it Proof of (a).} 
Take any spin isomorphism $\varphi$. 
By the Unlinked Section Lemma (b) of \S2 we can make a fiberwise isotopy so that 
$\varphi(\partial V_0)=\partial V_0'$ over $N_0$. 
Since $\pi_2(SO_2)=0$, we can then modify $\varphi$ over $N-N_0$ so that 
$\varphi(\partial V_0)=\partial V_0'$. 
\qed

\smallskip
{\it The Seifert Surface Lemma (b)} follows by (a) because 
$A_f=[\t f_*(V_0,\partial V_0)]$ by the Alexander Duality Lemma of \S2 (and the 
same for $f,V$ replaced by $f',V'$).

\smallskip
{\it Proof of (c).} 
Denote by $\widehat e\in H^2(N\times I,N\times\partial I)$ the obstruction to
extension of `the union' of normal sections of $\t f$ in $C_f$ and of $\t f'$ 
in $C_{f'}$ to a normal section of $Y$ in $M$.
Then $\widehat e$ goes to $\overline e$ under the composition 
$$H^2(N\times I,N\times\partial I)\overset{PD}\to
\cong H_2(N\times I)\cong H_2(N)\cong H_2(\partial V_0)\to H_2(Y).$$
Since the normal bundle of $N\cong\partial V_0$ in $Y$ is trivial, this 
implies that $\overline e\cap \overline e=0$.
\qed

\smallskip
(Clearly, $V_0\cong V$ and $V_0'\cong V'$, so $Y\cong V\cup_N(-V')$. 
Note that a spin isomorphism $\varphi$ such that 
$\varphi(\partial V_0)=\partial V_0'$ is unique by the Spin Framing Lemma of 
\S2 because $\pi_3(SO_2)=0$.)   

%(Although $\overline e\cap \overline e=0$ under the assumptions of the 
%Compression Theorem, this does not prove it because $\sigma(Y)$ could be odd 
%[Ta04, Ta06].) 

For an embedding $f:N\to\R^5$ denote by $if$ the composition of $f$ and the 
standard inclusion $\R^5\subset\R^6$.  
Recall that $d(W(if))=0$ because $2W(if)=0$.

\smallskip
{\bf Compressed Kreck Invariant Lemma.} 
{\it Let $f,f':N\to\R^5$ be embeddings such that $W(if)=W(if')$. 
For Seifert surfaces 
$$\t f:V\to\R^5\quad\text{and}\quad \t f':V'\to\R^5\quad\text{we have}\quad 
8\eta(if,if')=\sigma(V)-\sigma(V').$$}  

%\smallskip
{\it Proof.} 
Take a spin isomorphism $\varphi:\partial C_f\to\partial C_{f'}$ given by 
the Seifert Surface Lemma (a) and a submanifold $Y\subset M_\varphi$ 
given by the Seifert Surface Lemma (b). 
Then 
$$8\eta(if,if')=
\sigma(V\cup V')-\overline e\cap\overline e=\sigma(V)-\sigma(V')$$
by the Kreck Invariant Lemma, Novikov-Rokhlin additivity and the Seifert 
Surface Lemma (c).  
\qed

\bigskip
{\bf Proof of the Compression Theorems.}

\smallskip
{\bf Compression Theorem (an equivalent formulation).} 
{\it (1) For each $u\in H_1(N)$ a compressible embedding $f_0:N\to\R^6$ 
such that $W(f_0)=u$ exists if and only if $2u=0$. 

(2) Take an element $u\in H_1(N)$ of order 2 and a compressible embedding 
$f_0:N\to\R^6$ such that $W(f_0)=u$. 
A compressible embedding 
$$f:N\to\R^6\quad\text{such that}\quad W(f)=u\quad
\text{and}\quad \eta(f,f_0)=a$$ 
exists if and only if either $a$ is even or $u$ is not spin simple. }  

\smallskip
The above formulation is equivalent to the one in \S1 by the Classification 
Theorem.  

The new part of the Compression Theorem is the extension of (2) from the 
case when $N$ is a homology sphere to the general case (except for the 
`if' implication for $a$ even). 
\footnote{The `only if' part of (1) is easy, see [Vr89] or below; 
the `if' part of (1) follows by the Realization Lemma 
and the Relation Lemma (b) below, both are essentially proved in [SST02].
The case $N=S^3$ of (2) is known [Ha66, Ta06]; this case together with the 
Addendum to the Classification Theorem imply the `if' implication for $a$ even 
in (2), and an analogous statement holds for $(4k-1)$-manifolds in $\R^{6k}$.}   

%Let us explain in details which parts of the Compression Theorem are 
%known and which are new. 
%Thus the new (and non-trivial) parts of the Compression Theorem are 
%--- the `if' part in the non-spin-simple case of (a)  
%(or, equivalently, the non-spin-simple case of (2)); 
%--- generalization to $N$ not an integral homology sphere of the necessity of 
%(2) in the `only if' part in the spin simple case of (a) (or, equivalently, 
%such a generalization of the `only if' part of the spin simple case of (b2)). 
%A particular case is proved in [Ta06]: {\it an embedding 
%$f:N\to\R^6$ of an integral homology 3-sphere $N$ compresses to $\R^5$ if and 
%only if $\eta(f)\equiv\mu(N)\mod2$}. 
%Here $\eta(f)$ and $\mu(N)$ are the Haefliger and the Rokhlin invariants, 
%respectively. 

\smallskip
{\it Proof of the `only if' part in the Compression Theorem (1).}
The section $\xi$ of the normal bundle to $f$ that is normal to 
$\R^5\subset\R^6$, is unlinked.
The normal bundle of any codimension 2 embedding of an orientable manifold
is trivial.
Hence by the Bo\'echat-Haefliger Invariant Lemma of \S2 we have 
$\pm2W(f)=BH(f)=PDe(\xi^\perp)=0$.
\qed

\smallskip
{\it Proof of the Codimension Two Compression Theorem.} 
By the `only if' part of the Compression Theorem (1) $2W(f)=0=2W(f')$. 
So $W(f)=W(f')$. 

First we prove the case of $S^4$. 
The closure $V$ of a connected component of $S^4-f(N)$ is a Seifert surface of 
$f$. 
Define analogously $V'$. 
Since $V,V'\subset S^4\subset\R^5$, we have $\sigma(V)=\sigma(V')=0$. 
Hence $\eta(f,f')=0$ by the Compressed Kreck Invariant Lemma. 
 
The case of $S^2\times S^2$ is analogous for the following reasons:   

$\bullet$
any closed orientable 3-submanifold $U$ of $S^2\times S^2$ splits 
$S^2\times S^2$ 
(because in the opposite case $\#(U\cap\Sigma)=1$ for a certain circle 
$\Sigma\subset S^2\times S^2$ which is impossible), and 

$\bullet$
the signature of a non-closed connected oriented 4-submanifold $V$ 
of $S^2\times S^2$ is zero 
(because for the inclusion $i:V\to S^2\times S^2$ 
we have $x\cap y=i_*x\cap i_*y$, so there is a basis $\{x_1,\dots,x_s\}$ in 
$H_2(V;\Q)$ such that all the intersections $x_i\cap x_j$ except possibly 
$x_{s-1}\cap x_s$ are zeroes).  
\qed

\smallskip
Now we complete the proof of the Compression Theorem in its alternative 
formulation. 

In the rest of this section let $f:N\to\R^5$ be an embedding. 
Denote by 

$\bullet$
$C_f$ the closure of the complement in $S^5\supset\R^5$ to a tubular
neighborhood of $f(N)$ and   

$\bullet$
$\nu_f:\partial C_f\to N$ the restriction of the normal bundle of $f$.  

By a {\it section} we always mean a section $N\to\partial C_f$ of 
$\nu_f:\partial C_f\to N$.

%\smallskip
%{\it Definition of the spin structure $s_f$ on $N$ corresponding to an 
%embedding $f:N\to\R^5$.} 
Take the normal section of $f$ pointing to $\t f(V)$, where $\t f:V\to\R^5$
is an extension of $f$ to an embedding of a connected oriented 4-manifold $V$ 
with boundary $\partial V=N$; such a {\it Seifert surface} exists by [Ki89, 
VIII, Theorem 3]. 
\footnote{In other words, take the normal section $\xi:N\to\R^5-f(N)$ of $f$ 
which is {\it unlinked}, i.e. $\xi_*[N]=0\in H_3(\R^5-f(N))$; the existence and 
uniqueness of unlinked section follows by the analogue of Difference Lemma (a) 
of \S2 because for fixed section $\xi$ the map $\zeta\mapsto d(\xi,\zeta)$ 
defines a 1--1 correspondence between sections and $H_2(N)$; the existence of 
an unlinked section is the first step in the proof of the existence of a 
Seifert surface.} 
Complete this section to a normal framing on 
$f(N)$ in $\R^5$ so that the fixed orientation on $N$ followed by the orientation of the normal bundle 
(obtained from the framing) forms the fixed orientation on $\R^5$. 
For an embedding $f:N\to\R^5$ let $s_f$ be {\it the spin structure on $N$ 
defined by the composition of the constructed (`unlinked') normal framing and 
the standard framing of $\R^5\subset\R^8$}.  

%A section $\xi:N\to\partial C_f$ for $f$ is called {\it unlinked} if 
%$in_*\xi_*[N]=0\in H_3(C_f)$. 
%(Cf. the definition of an unlinked section of an embedding $f:N\to\R^6$ in \S2.)
%\smallskip
%{\bf Unlinked Section Lemma.}
%{\it An unlinked section exists and is unique up to fiberwise homotopy.} 
%%If $\xi$ is an unlinked section, then $\xi_*[N]=\partial A_f$.}  

%\smallskip
Fix the spin structure $s_g$ on the manifold $N$, where $g$ is the embedding 
used in the definition of the Whitney invariant. 
Then {\it spin structures are identified with elements of $H_2(N;\Z_2)$} and we 
write $s_f\in H_2(N;\Z_2)$.  

\smallskip
{\bf Relation Lemma.} {\it Let $f,f':N\to\R^5$ be embeddings. 

(a) If $W(if)=W(if')$, then $\eta(if,if')\equiv\mu(s_f)-\mu(s_{f'})\mod2$. 

(b) $W(if)=\beta s_f$, where $\beta:H_2(N;\Z_2)\to H_1(N)$ is the (homology) 
Bockstein homomorphism.} 

\smallskip
{\it Proof of (a).} 
Take an embedding $f:N\to\R^5$ and its Seifert surface $\t f:V\to\R^5$. 
The normal section of $\t f$ compatible with the orientations defines a spin 
structure $s_{\t f}$ on $V$. 
We have $s_{\t f}|_N=s_f$. 
Therefore the required congruence holds by the Compressed Kreck Invariant 
Lemma. 
\qed

\smallskip
{\it Proof of (b).} 
Define $c(f)\in H_1(N)$ to be the homology class of the singular 1-submanifold 
of a general position homotopy $H$ between $g|_{N_0}$ and $f$ (this is a 
submanifold by general position; it is well-known that $c(f)$ is indeed 
independent of $H$). 
This definition agrees with [SST02, Definition 3.3(1)] by [SST02, Remark 3.2]. 
We have   
$$W(if)\underset{(1)}\to=c(f)\underset{(2)}\to=\beta\varphi_{f,\nu,*}y
\underset{(3)}\to=\beta s_f.$$ 
Here (1) is the commutativity of [Vr89, Theorem 3.1] which holds for $k=0$, 
cf. [Ya83],  

(2) is [SST02, Lemma 3.5], cf. [Sk05, the Whitney Invariant Lemma ($\beta$)] 
and 

(3) is clear for the fixed spin structure $s_g$ on $N$, the unlinked framing 
$\nu$ of $f$ and $\t c_f=\varphi_{f,\nu,*}y=s_f\in H_2(N;\Z_2)$ defined 
in [SST02, after Definition 3.3]. 
\qed

\smallskip
{\it Proof of the `only if' part of the Compression Theorem (2).} 
Suppose that $u$ is spin simple and let us prove that $a$ is even. 
Take embeddings $F,F_0:N\to\R^5$ such that $f=iF$ and $f_0=iF_0$ ($F$ and $F_0$ 
need not be unique up to isotopy). 
Since $W(iF)=W(iF_0)=u$, by the Relation Lemma (b) 
$\beta(s_F)=\beta(s_{F_0})=u$. 
So $a=\eta(iF,iF_0)\equiv\mu(s_F)-\mu(s_{F_0})=0\mod2$ 
(by the Relation Lemma (a) and because $u$ is spin simple).  
\qed

\smallskip
{\bf Realization Lemma.} 
{\it For each $x\in H_2(N;\Z_2)$ there is an embedding $f:N\to\R^5$ such that 
$s_f=x$. }

\smallskip
{\it Proof.} The proof is analogous to that of [SST02, Theorem 3.8] but we 
present the details because the Realization Lemma does not follow from the 
{\it statement} of [SST02, Theorem 3.8].  

By [Ka79] there exists a 4-manifold $V$ and a spin structure $s$ on $V$ 
such that $\partial V=N$, $s|_N=x$ and $V$ has a handle 
decomposition that consists of one 0-handle and some 2-handles attached 
to the 0-handle simultaneously. 
Since $V$ is a non-closed connected spin 4-manifold, it is parallelizable. 
So $V$ immerses into $\R^5$. 
Furthemore, since $V$ has a 2-dimensional spine, the immersion is regular 
homotopic to an embedding. 
Let $f$ be the restriction to $N$ of such an embedding. 
Now the Lemma follows because $s_f=s|_N=x$. 
\qed

\smallskip
{\it The 'if' part of the Compression Theorem (1)} follows from the Relation 
Lemma (b) and the Realization Lemma. 

\smallskip
{\bf Double Trefoil Knot Lemma.} 
{\it The knot $t\#t:S^3\to\R^6$ is compressible} [Ha66]. 

\smallskip
{\it Proof.} 
Let $V$ be a closed simply-connected spin 4-manifold such that $\sigma(V)=16$ 
and $V$ has a handle decomposition that consists of one 0-handle and some 
2-handles attached to the 0-handle simultaneously (e.g. the Kummer surface). 
Analogously to the proof of the Realization Lemma there is an embedding 
$f:V_0\to\R^5$. 
Then $\eta(if|_{\partial V_0})=\sigma(V_0)/8=2$ by the Compressed 
Kreck Invariant Lemma (because $\eta(f')=\sigma(B^4)/8=0$ for the standard 
embedding $f'$). 
Hence $t\#t$ is isotopic to $if|_{\partial V_0}$ and so is compressible. 
\qed 

\smallskip
{\it Proof of the `if' part of the Compression Theorem (2).} 
The case when $a$ is even follows by taking $f:=f_0\#at$ which is compressible 
by the Double Trefoil Knot Lemma. 
So suppose that $u$ is not spin simple. 
Then there exist spin structures $s,s'$ on $N$ such that 
$\beta(s)=\beta(s')=u$ and $\mu(s)\ne\mu(s')$. 
By the Realization Lemma there are embeddings $f,f':N\to\R^5$ such that 
$s_f=s$ and $s_{f'}=s'$ ($f$ is not necessarily the one required in the 
Theorem).
By the Relation Lemma (b) and (a) 
$$W(if)=\beta(s)=u=\beta(s')=W(if')\quad\text{and}\quad 
\eta(if,if')\equiv\mu(s)-\mu(s')=1\mod2.$$ 
Hence by the Classification Theorem and Addendum 
$if=(if')\#lt$ for some odd $l$. 
Hence 
$$W^{-1}(u)=\{(if)\#kt\ |\ k\in\Z\}=\{(if)\#2kt,(if')\#2kt\ |\ k\in\Z\}.$$ 
This and the Double Trefoil Knot Lemma imply that 
$W^{-1}(u)$ consists of compressible embeddings. 
Hence we can take $f_0\#at$ as the required compressible embedding.

%\newpage
\head 4. A smoothing proof of the Classification Theorems \endhead

For an alternative proof of the Classification Theorem we give an alternative 
definition of the Kreck invariant $\eta^H$ and proof of the Injectivity Lemma 
with $\eta$ replaced by $\eta^H$ (the Classification Theorem follows from 
the Injectivity Lemma with $\eta$ replaced by $\eta^H$ in the same way as in 
\S2). 

Clearly, every smooth (i.e.  differentiable) map is piecewise differentiable.
The forgetful map from the set of piecewise linear embeddings (immersions)
up to piecewise linear isotopy (regular homotopy) to the set of piecewise
differentiable embeddings (immersions) up to piecewise differentiable isotopy
(regular homotopy) is a 1--1 correspondence [Ha67].
Therefore we can consider any smooth map as PL one, although this is
incorrect literally.

Analogously to \S1 we define the PL Whitney invariant 
$W_{PL}:\Emb\phantom{}^6_{PL}(N)\to H_1(N)$ from the set of PL isotopy classes 
of PL embeddings $N\to\R^6$. 

\smallskip
{\bf PL Classification Theorem.}
{\it $W_{PL}$ is bijective} [Vr77, Theorem 1.1], cf. 
[Hu69, \S11, BH70, Th\'eor\`eme 1.6, Sk97, Corollary 1.10, 
Sk07, Theorem 2.8.b], \S5, (5).

\smallskip
{\it An alternative definition of the Kreck invariant.}  
Let $f,f':N\to\R^6$ be two (smooth) embeddings such that $W(f)=W(f')$.
Then by the PL Classification Theorem there is a PL isotopy 
$F:\R^6\times I\to\R^6\times I$ between $f$ and $f'$.
Making a PL isotopy of $\R^6\times I$ we may assume that $F$ smooth outside a 
fixed single point [Ha67], cf. beginning of [BH70, Bo71].
Consider a small smooth oriented 7-ball with the center at the image of this 
point. 
Let $\Sigma$ be the boundary of this ball.
Take the natural orientation on the 3-sphere $F^{-1}\Sigma$.
Denote by $\overline F:F^{-1}\Sigma\to\Sigma$ the abbreviation of $F$. 
Define $\eta(\overline F)\in\Z$ as in the beginning of \S2. 
Define the Kreck invariant by 
$\eta^H_{W(f)}(f,f'):=\eta(\overline F)\mod d(W(f))$. 
It is well-defined by the Smoothing Lemma (1) below (instead of the 
Independence Lemma).

%is independent of the PL isotopy of $F$ and is the complete obstruction to the 
%existence of a PL isotopy, relative to the boundary, from $F$ to a smooth 

\smallskip
A PL {\it concordance} between two smooth embeddings $f,f':N\to\R^6$ is a PL
embedding $F:N\times I\to\R^6\times I$ such that $F(x,0)=(f(x),0)$ and
$F(x,1)=(f'(x),1)$ for each $x\in N$.

\smallskip
{\bf Smoothing Lemma.}
{\it Let $F$ be a PL concordance between two smooth embeddings $f$ and $f'$.
Then

(1) for each PL concordance $F'$ between $f$ and $f'$ the integer
$\eta(\overline F)-\eta(\overline F')$ is divisible by $d(W(f))$, and

(2) for each $s\in\Z$ there exists a PL concordance $F'$ between
$f$ and $f'$ such that $\eta(\overline F)-\eta(\overline F')=sd(W(f))$.}

\smallskip
The Smoothing Lemma is proved below using the Bo\'echat-Haefliger formula for 
$\eta(\overline F)$. 
%the smoothing obstruction.

\smallskip
{\it Proof of the Addendum with $\eta$ replaced by $\eta^H$.} 
From the cone on the Haefliger trefoil knot $t:S^3\to\R^6$ we obtain a PL 
isotopy $T$ between $t$ and the standard embeding such that 
$\eta(\overline T)=1$.  
From the identical isotopy $I$ between $f$ and $f$ we obtain an isotopy 
$I\#kT$ between $f$ and $f\#kt$.  
Hence 
$$\eta^H(f\#kt,f)=\eta(\overline{I\#kT})=\eta(\overline{kT})\equiv k
\mod d(W(f)).\quad\qed$$

%\smallskip
{\it Proof of the equivalence of two definitions of the Kreck invariant.}
If $f,f':N\to\R^6$ are two (smooth) embeddings such that $W(f)=W(f')$, then by 
the Classification Theorem and the Addendum $f'=f\#kt$ for some 
$k=\eta(f',f)\mod d(W(f))$. 
Now $\eta(f',f)=k=\eta^H(f\#kt,f)=\eta^H(f',f)$ by the $\eta$- and 
$\eta^H$-versions of the Addendum. 
\qed 

\smallskip
{\it An alternative proof of the Injectivity Lemma with $\eta$ replaced by 
$\eta^H$.} 
Since $W(f)=W(f')$, by the PL Classification Theorem there
is a PL isotopy $F$ between $f$ and $f'$.
Since $\eta^H(f,f')\equiv0\mod d(W(f))$, by the Smoothing Lemma (2) 
there exists a PL concordance $F_1$ between $f$ and $f'$ such that 
$\eta(\overline F_1)=0$.
Then $F_1$ is PL concordant relative to the boundary to a smooth isotopy.
Thus $f$ is isotopic to $f'$. \qed

\smallskip
{\it Proof of (1) in the Smoothing Lemma.}
Two PL isotopies $F$ and $F'$ between $f$ and $f'$ together form an embedding
$\Psi:=F\cup F':N\times S^1\to\R^7$.
Take a ball $B^4\subset Q:=N\times S^1$. 
Set $Q_0:=Q-\Int B^4$. 
Denote by 

$\bullet$
$C_\Psi$ the closure of the complement in $S^7$ to a tubular neighborhood of 
$\Psi(Q)$;  

$\bullet$
$|\cdot,\cdot|$ the distance in $Q$ such that $B^4$ is a ball of radius 2. 

For a section $\Xi:Q_0\to\partial C$ define the map 
$$\overline\Xi:Q\to S^7-\Psi(Q_0)\quad\text{by}\quad
\overline\Xi(x)=\cases \xi(x)&x\in Q_0\\ 
\Psi(x)&|x,Q_0|\ge1\\ 
|x,Q_0|\Psi(x)+(1-|x,Q_0|)\xi(x)&|x,Q_0|\le1\endcases.$$
A section $\Xi:Q_0\to\partial C_\Psi$ is called {\it unlinked} if 
$\overline\Xi_*[Q]=0\in H_4(S^7-\Psi Q_0)$, cf. \S2, [BH70, KS05]. 
By [HH63, 4.3, BH70, Proposition 1.3, Lemme 1.7] 
{\it an unlinked section exists and an unlinked section is unique on 
the 2-skeleton of $Q$.}
%ne ispolzuetsja
%We have $PDe(\xi^\perp)\mod2=PDw_2(M)$ (which could now be non-zero).

The isotopies $F$ and $F'$ are ambient and we may assume that the corresponding
enveloping isotopies $\R^6\times I\times I\to\R^6\times I\times I$ are smooth
outside some balls in $\R^6\times[\frac13,\frac23]\times[\frac13,\frac23]$.
Hence an unlinked section of $f$ can be extended to an unlinked section $\Xi$ 
of $F\cup F'$.
%!!!
Identify $H_2(N\times S^1)$ with $H_1(N)\oplus H_2(N)$ by the K\"unneth 
isomorphism.
Now (1) follows because
$$\eta(\overline F)-\eta(\overline F')\underset{1}\to=\eta(\overline{F\cup F'})
\underset{2}\to=\frac{PDe(\Xi^\perp)^2}8-\frac{p_1(N\times S^1)}{24}
\underset{3}\to=\frac{(2W(f)\oplus2\alpha)^2}8\underset{4}\to=
W(f)\bigcap\limits_N\alpha$$
for some $\alpha\in H_2(N)$.
Here the first equality is clear.
The second follows by [BH70, Th\'eor\`eme 2.1, cf. Bo71, Fu94] because the 
Kreck invariant coincides with the Haefliger invariant for $N=S^3$ by [Wa66] or 
[GM86, Remarks to the four articles of Rokhlin, II.2.7 and III.excercises.IV.3, 
Ta04]. 

Let us prove the third equality.
Since $N\times S^1$ is parallelizable, it follows that $p_1(N\times S^1)=0$.
Thus it suffices to prove that
$PDe(\Xi^\perp)=\pm2(W(f)\oplus\alpha)$.
Since $e(\Xi^\perp)\mod2=w_2(N\times S^1)=0$, it follows that the projection 
of $PDe(\Xi^\perp)$ onto the second summand $H_2(N)$ is indeed an even element. 
The projection of $e(\Xi^\perp)$ onto the first summand is
$$PDe(\Xi^\perp)\cap[N\times1]=PDe(\Xi^\perp|_{N\times1})=BH(f)=\pm2W(f)$$
by the Bo\'echat-Haefliger Invariant Lemma of \S2. 

Let us prove the fourth equality.
We can represent elements
$$W(f)\oplus 0,\ 0\oplus\alpha\in H_2(N\times S^1) 
\quad\text{by}\quad Z\times S^1\quad\text{and}\quad L\times1,$$
respectively, where $Z$ is an oriented circle in $N$ and $L$ is an oriented 
sphere with handles in $N$.
We have $(0\oplus\alpha)^2=[L\times1]\cap[L\cap i]=0$.
Clearly, there is a circle $Z'$ in $N$ homologous to $Z$ and disjoint with $Z$.
Thus $(W(f)\oplus0)^2=[Z\times S^1]\cap[Z'\times S^1]=0$.
We also have $[L\times1]\cap_{N\times S^1}[Z\times S^1]=[L]\cap_N[Z]$.
All this implies the fourth equality.
\qed

\smallskip
{\it Proof of (2) in the Smoothing Lemma.}
Let $B^4\subset N\times (0,1)$ be a ball and denote
$C'=S^6\times I-F(N\times I-\Int B^4)$.
By general position, $C'$ is simply-connected.
By Alexander and Poincar\'e duality
$$H_i(C')\cong H^{6-i}(N\times I-\Int B^4,N\times \{0,1\})
\cong H_{i-2}(N\times I,\partial B^4)\cong H_{i-2}(N).$$
Thus $C'$ is 2-connected.
Hence the Hurewicz homomorphism $\pi_4(C')\to H_4(C')\cong H_2(N)$ is an 
epimorphism.
Therefore analogously to the construction at the end of \S1 (or to [BH70, proof 
of Theorem 1.6], see 
%[Sk06, 
\S5, (5)) and using Mayer-Vietoris sequence, for 
each $\alpha\in H_2(N)$ we can construct $F'$ and an unlinked section of $\Xi$ 
of $F$ extending an unlinked section of $f$ so that 
$PDe(\Xi^\perp)=BH(f)\oplus2\alpha$.
%!!!
Then $\eta(\overline F)-\eta(\overline F')=\pm W(f)\cap\alpha$ by the proof of
the Smoothing Lemma (1).
Now (2) follows because by Poincar\'e duality $W(f)\cap H_2(N)=d(W(f))\Z$.
\qed

\smallskip
{\it Proof of the Higher-dimensional Classification Theorem (b).}
The proof is analogous to the proof of the Classification Theorem. 
We use $\Emb^{6k}(S^{4k-1})\cong\Z$ [Ha62, Ha66] together with  the
higher-dimensional analogues of the PL Classification Theorem 
[Hu69, \S12, Vr77], the Bo\'echat-Haefliger Invariant Lemma and the Smoothing 
Lemma.
The latter result is proved analogously using [Bo71, Th\'eor\`eme 5.1]
instead of [BH70, Th\'eor\`eme 2.1].
We replace $p_1$ by $p_k$ and $w_2$ by $\overline w_{2k}$. 
Since 
$$N\times S^1=\partial(N\times D^2),\quad\text{we have}
\quad p_k(N\times S^1)=0\in H^{4k}(N\times S^1)\cong\Z.$$ 
We have $\overline w_{2k}(N\times S^1)=\overline w_{2k}(N)=0$ by the product 
formula for Stiefel-Whitney classes and by the following Lemma. 
\qed

\smallskip
{\bf Lemma.} {\it If $N$ is a closed $\Z_2$-homologically $(2k-2)$-connected 
$(4k-1)$-manifold, then all the mod 2 Stiefel-Whitney classes of $N$ (both 
tangential and normal) are zeros.} 

\smallskip
{\it Proof} (communicated to the author by D. Crowley). 
We prove the lemma for the tangential classes, which by the Whitney-Wu formula 
implies the case of normal classes. 
The only non-trivial cohomology groups $H^i(N;\Z_2)$ are those for $i=2k-1$ 
and $i=2k$. 
So we only need to prove that $w_{2k-1}(N)=0$ and $w_{2k}(N)=0$. 

Since $N$ is $(2k-2)$-connected, by the Wu formula [MS74, Theorem 11.14] 
we have $0=\Sq^1w_{2k-2}=w_{2k-1}(N)=v_{2k-1}(N)$. 
For the same reason and since $2\cdot2k>4k-1$, we have 
$w_{2k}(N)=v_{2k}(N)+\Sq^1v_{2k-1}(N)=0$. 
(Here $v_i(N)$ are the Wu classes.)
\qed

\smallskip
{\it Proof of the Higher-dimensional Classification Theorem (a).}
For 
$$1\le p\le2k-2\quad\text{we have}\quad 2\cdot6k\ge3p+2(4k-1-p)+4\quad\text
{and}\quad 6k\ge2p+4k-1+3.$$ 
Hence by [Sk02, Theorem 1.3.$\alpha$, Sk07, Group Structure Theorem 3.7, 
Sk06, Group Structure Theorem 2.1] we have the exact sequence of {\it groups}
$$0\to\Z\overset\zeta\to\to\Emb\phantom{}^{6k}(S^p\times S^{4k-1-p})
\overset\alpha\to\to \pi_{4k-1-p}(V_{2k+p+1,p+1})\to0.$$
The map $\zeta$ defines the action of $\Emb^{6k}(S^{4k-1})$ by embedded 
connected summation. 
It is injective because analogously to the proof of the Smoothing Lemma (1) we 
construct $F,F',\Xi$ and obtain 
$$\eta(\overline F)-\eta(\overline F')\underset{1}\to=\eta(\overline{F\cup F'})
\underset{2}\to=
\frac{PDe(\Xi^\perp)^2}8-\frac{p_k(S^p\times S^{4k-1-p}\times S^1)}{24}
\underset{3}\to=0.$$
Here the second equality holds by [Bo71, Th\'eor\`eme 5.1] while the third
equality holds because
$PDe(\Xi^\perp)\in H_{2k}(S^p\times S^{4k-1-p}\times S^1)=0$ and because
$S^p\times S^{4k-1-p}\times S^1$ is parallelizable (since either $p$ or
$4k-1-p$ are odd).

This sequence splits because $\alpha$ has a right inverse $\tau$
[Sk02, Torus Lemma 6.1, Sk07, Torus Lemma 6.1].
\qed

%\newpage
\head 5. Some remarks and conjectures \endhead

\smallskip 
{\bf Constructions of embeddings.}

\smallskip 
(1) The following idea could
perhaps be used to construct a map $[N;S^2]\to\Emb^6(N)$. Given a map
$\varphi:N\to S^2$ take the map $\varphi\times\pr_2:N\times S^2\to
S^2=\C P^1\subset\C P^\infty$, prove that the latter is
zero-bordant, use the zero-bordism as a model for $C_f$ and
analogously to [Fu94] by surgery obtain $C_f$ and $f$. 

\smallskip
(2) {\it The Hopf construction of an embedding $\R P^3\to S^5$.}
Represent $\R P^3=\{(x,y)\in\C^2\ |\ |x|^2+|y|^2=1\}/\pm1$.
Define $f:\R P^3\to S^5\subset\C^3$ by $f[(x,y)]=(x^2,2xy,y^2)$.
It is easy to check that $f$ is an embedding.
(The image of this embedding is given by the equations $b^2=4ac$, 
$|a|^2+|b|^2+|c|^2=1$.)

It would be interesting to obtain an explicit construction of an embedding 
$\R P^3\to S^5$ whose composition with the standard inclusion $S^5\subset\R^6$
is not isotopic to the Hopf embeding. 

\smallskip
(3) {\it An embedding $\C P^2\to\R^7$} [BH70, p. 164]. 
It suffices to construct {\it an embedding $f_0:\C P^2_0\to S^6$ such that the 
boundary 3-sphere is the standard one.}  
Recall that $\C P^2_0$ is the mapping cylinder of the Hopf map $h:S^3\to S^2$. 
Recall that $S^6=S^2*S^3$. 
Define $f_0[(x,t)]:=[(x,h(x),t)]$, where $x\in S^3$. 
In other words, the segment joining $x\in S^3$ and $h(x)\in S^2$ is mapped onto 
the arc in $S^6$ joining $x$ to $h(x)$.   

By [Ta06, Proposition 3.7], cf. the above construction, there is a 
Seifert surface $\t g':\C P^2_0\to\R^6$ of the standard embedding $S^3\to\R^6$ 
such that for the standard Seifert surface $\t g:D^4\to\R^6$ we have 
$\overline e=[\C P^1]\in H_2(\C P^2_0)$ and $\overline p_1=1$.  
Hence for each two embeddings $f,f_0:N\to\R^6$ such that $W(f)=W(f_0)$

--- there are Seifert surfaces $\t f:V\to\R^6$ and $\t f_0:V_0\to\R^6$ such 
that for $\varphi$ as in the Seifert Surface Lemma (a) we have $\overline p_1=0$, 
so that $\eta(f,f_0)\equiv\sigma(Y)/8\mod d(W(f)).$

--- there are Seifert surfaces $\t f:V\to\R^6$ and $\t f_0:V_0\to\R^6$ such 
that for $\varphi$ as in the Seifert Surface Lemma (a) $\sigma(Y)=0$, so that 
$\eta(f,f_0)\equiv-\overline p_1/8\mod d(W(f))$. 

\smallskip
(4) {\it An explicit construction of the generator $t\in\Emb^6(S^3)$} 
[Ha62, 4.1].
Denote coordinates in $\R^6$ by $(x,y,z)=(x_1,x_2,y_1,y_2,z_1,z_2)$.
The higher-dimensional trefoil knot $t:S^3\to\R^6$ is obtained by joining
with two tubes {\it the higher-dimensional Borromean rings}, i.e. the three
spheres given by the following three systems of equations:
$$\cases x=0\\ |y|^2+2|z|^2=1\endcases, \qquad
\cases y=0\\ |z|^2+2|x|^2=1\endcases \qquad\text{and}\qquad \cases
z=0\\ |x|^2+2|y|^2=1 \endcases.$$

Let us sketch the proof of the surjectivity of $\eta:\Emb^6(S^3)\to\Z$, cf. 
[Ta04]. 
We use the definition of $M$ and $y$ from the beginning of \S2.
It suffices to prove that $\eta(t)=1$ for the higher-dimensional trefoil knot
$t$.
The Borromean rings $S^1\sqcup S^1\sqcup S^1\to\R^3$ span three 2-disks
whose triple intersection is exactly one point.
Analogously, the higher-dimensional Borromean rings
$S^3\sqcup S^3\sqcup S^3\to\R^6$ span three 4-disks whose triple
intersection is exactly one point, cf. [Ha62, Ta06, Sk06, Sk07].
Take the higher-dimensional trefoil knot $t:S^3\to\R^6$ obtained by joining
Borromean rings with two tubes.
For its framing take the section formed by the first vectors of the framing.
This section spans an immersed 4-disk in $C_t$ whose triple self-intersection
is exactly one point.
The union of this disk with the disk
$$x\times B^4\subset\partial D^3\times B^4\quad\subset\quad
(S^6-tS^3\times\Int D^3)
\bigcup\limits_{tS^3\times \partial D^3=\partial B^4\times \partial D^3}
B^4\times\partial D^3=M$$
is an immersed 4-sphere representing the class $y\in H_4(M)$.
Recall that we take a framing $\varphi$ of $f$ so that $yPDp_1(M)=0$.
Then $6\eta(t)=y^3=6$. 
\qed

\smallskip
(5) {\it An alternative proof of the surjectivity of $W$} [Vr77]
This proof, although more complicated, is interesting because unlike the 
previous one it can be generalized to the proof of the injectivity of $W_{PL}$. 

%By the Whitney Embedding Theorem there exists an embedding $f:N\to\R^6$.
Let $C'_g$ be the closure of the complement in
$S^6$ to the space of the normal bundle to $N$ restricted to $N_0$
(i.e. to the regular neighborhood of $gN_0$ modulo $g\partial B^3$).
By Alexander and Poincar\'e duality $H_3(C'_g)\cong H^2(N_0)\cong H^2(N)\cong H_1(N)$.
By general position and Alexander duality $C'_g$ is 2-connected.
Hence the Hurewicz homomorphism $\pi_3(C'_g)\to H_3(C'_g)$ is an isomorphism.

So for each $u\in H_1(N)$ there is a map $\bar u:S^3\to C'_g$ whose homotopy
class goes to $u$ under the composition of the above isomorphisms.
Let $\hat u:B^3\to C_g'$ be a connected sum of $g|_{B^3}:B^3\to C_g'$ and $\bar u$.
Since $C'_g$ is simply-connected, we can modify $\hat u$ by a homotopy relative to
the boundary to an embedding $f':B^3\to C_g'$ using Whitney trick.
Define $f$ to be $g$ on $N_0$ and $f'$ on $B^3$.
(Note that $f$ depends on $f'$ whose isotopy class is not unique.)

Clearly, $f$ is smooth outside $\partial B^3$.
By a slight modification of $f'$ we may assume that $f$ is smooth on
$\partial B^3$ (see the details in [HH63, bottom of p. 134]).
Clearly, the class of $f|_{B^3}\cup g|_{B^3}$ in $H_3(C'_g)$ goes to $u$ under
$H_3(C'_g)\cong H^2(N_0)\cong H^2(N)\cong H_1(N)$.
Hence $W(f)=u$.
\qed

\smallskip
It would be interesting to construct {\it a map} $\theta:H_1(N)\to\Emb^6(N)$ 
such that $W\theta(u)=u$, and thus an {\it absolute} Kreck invariant. 
For this we would need to make the above construction or the construction 
at the end of \S1 in a {\it canonical} way (i.e. choose canonical embedding 
$u$, extension $\overline u:D^2\to S^6$ and the 
2-framing normal to $\overline u(D^2)$).  

\smallskip
{\it Proof of the PL Classification Theorem.}
The surjectivity of $W_{PL}$ is proved as above (using the
Penrose-Whitehead-Zeeman Embedding Theorem [Hu69, Sk07, \S2] instead of the
Whitney trick).

The injectivity of $W_{PL}$ follows because in the proof of the surjectivity
of $W$  the embedding $f':B^3\to C_f'$ is unique up to homotopy, so the
embedding $f'$ is unique up to PL isotopy by Zeeman Unknotting Theorem [Hu69]
because $C_g'$ is 2-connected. \qed

\smallskip
I am grateful to Jacques Bo\'echat for indicating that the
injectivity in [Bo71, Theorem 4.2] is wrong without the assumption
that $H_{k+1}(N)$ has no 2-torsion.
%(which assumption though holds for some $n$ and $k$ automatically)?
This theorem
%classifies embeddings of closed $k$-connected
%orientable $n$-manifolds into $S^{2n-k}$.
states that {\it for $n-k$ odd $\ge3$, is a closed orientable
$k$-connected $n$-manifold $N$ the Bo\'echat-Haefliger invariant
is a 1--1 correspondence between the set of PL isotopy classes of
PL embeddings $N\to\R^{2n-k}$ and $\rho_2^{-1}\bar w^{n-k-1}(N)$,
where $\rho_2$ is the reduction modulo 2.}
%(Of course in this result the PL category can be replaced by smooth for
%$n\ge2k+4$.)
The assumption was used in the proof of the injectivity on p. 150,
line 3 from the bottom, in order to conclude that $0 = \chi_\sigma
- \chi_{\sigma'} = \pm 2 d(\sigma,\sigma')$ implies that
$d(\sigma,\sigma')=0$.
%If we take $n=3$, $k=0$ and $M=RP^3$, then the above Theorem 4.2 states that
%every two PL embeddings $\R P^3\to S^6$ are PL isotopic.
%We can analogously construct a higher-dimensional examples.
%These contradict to [HH63, Theorem 2.4, Hu69, Chapter 11, Ba75, Theorem, Vr77,
%Theorem 1.1] (note that the PL cases of these results follow from [Sk97, Sk02]).

\bigskip
{\bf Compression problem.}

\smallskip
{\it A direct proof that if $f(N)\subset\R^5$ and $W(f)=W(f_0)$, then 
$\eta(f,f_0)$ is even, for $N$ being a connected sum of some copies of 
$S^1\times S^2$ and a homology 3-sphere.} 
By [Ki89, VIII, Theorem 3] there are Seifert surfaces $\t f:V\to\R^5$ of $f$ 
and $\t f_0:V_0\to\R^5$ of $f_0$. 
Take $\varphi$ and $Y$ given by the Seifert Surface Lemma (a) and (b). 
We have 
$$0=w_2(M)|_Y=w_2(Y)+w_2(Y\subset M)=w_2(Y)+\rho PD\overline e,
\quad\text{so}\quad\rho\overline e=PDw_2(Y).$$ 
Since $V$ has codimension 1 in $\R^5$, it follows that $\t f$ has a normal 
section.
Analogously $\t f_0$ has one.
Since for our manifold $N$ {\it every class in $H_2(N)$ is realizable by an 
embedding of a disjoint union of spheres $S^2$}, by the Seifert Surface Lemma 
(c) the class $\overline e$ is realizable by an embedding 
$S^2\to N\cong\partial V_0\subset Y$. 
Now the Theorem follows by the Kreck Invariant Lemma and the corollary of the 
Rokhlin Theorem [Ma80, Corollary 1.13]. 
\qed

\smallskip
We conjecture that {\it any two embeddings $N\to\R^4$ of a closed 
2-manifold $N$ whose images are in $\R^3\subset\R^4$ are isotopic} (i.e. that 
only the standard embedding $N\to\R^4$ is compressible). 
It would be interesting to check whether the Hudson torus 
$S^1\times S^1\to\R^4$ [Sk07, \S3] is compressible. 

We conjecture that {\it any two embeddings $N\to\R^6$ whose images are in
$S^2\times S^2\subset\R^6$ are isotopic}, i.e that the Codimension Two 
Compression Theorem holds without the assumptions. 

It would be interesting to construct an example of a 3-submanifold of $S^4$ 
with non-trivial 2-torsion in homology. 
Such a torsion is necessarily of the form $G\oplus G$ for some group $G$ (so 
$\R P^3$ does not embed into $S^4$). 
Moreover, the restriction of linking form onto each $G$-summand is trivial 
(so $\R P^3\#\R P^3$ does not embed into $S^4$). 
It would be interesting to know which forms $T\times T\to\Q/\Z$ are realizable 
as linking forms of 3-manifolds.

%One of the systems of equations defining a 3-manifold appearing in the theory 
%of integrable systems is 
%$$R_1^2+R_2^2+R_3^2=c_1,\quad R_1S_1+R_2S_2+R_3S_3=c_2,\quad
%\frac{S_1^2}{A_1}+\frac{S_2^2}{A_2}+\frac{S_3^2}{A_3}=c_3,$$ where
%$R_i$ and $S_i$ are variables while $A_1<A_2<A_3$ and $c_i$ are
%constants.
%This system of equations corresponds to the Euler integrability case, see other 
%systems in [BF04, Chapter 14].

It would be interesting to characterize spin simple classes in 3-manifolds. 

The PL analogue of the Compression Theorem (in the formulation of \S1) for the 
non-spin-simple case 

is true (even without the non-spin-simplicity assumption), 

follows by the Compression Theorem (1) and the PL Classification Theorem 
of \S4,  

thus is essentially known. 

The analogue of the Compression Theorem (1) for the PL category and 
$(4k-1)$-manifolds in $\R^{6k}$ holds by [Vr89] (even without the 
non-spin-simplicity assumption).

{\it If $N$ is a 3-manifold with non-empty boundary such that $H_1(N)$ has no 
torsion and $f,f':N\to\R^5$ are two embeddings whose images are contained in 
$\R^4$, then $f$ and $f'$ are isotopic} by [Sa99, Theorem 3.1]. 
(Indeed, take the natural trivialization $\tau$ of the normal bundles of $f$ 
and $f'$ whose first vectors are orthogonal to $\R^4\subset\R^5$. 
Define the map $i^\tau:N\to\R^5$ as the shift of $f$ by vector orthogonal to 
$\R^4\subset\R^5$. 
Define the {\it Seifert form} $L^{\tau,f}:H_2(N)\times H_2(N)\to\Z$ by 
$L^{\tau,f}(\alpha,\beta):=\link(f_*\alpha,i^\tau_*\beta)$. 
Clearly, $L^{\tau,f}=0$ and the same for $f'$.)  
%This again implies that $W(f)=W(f')$ by the Compression Theorem (1). 

It would be interesting to know whether {\it an embedding $f:N_0\to S^5$ 
extends to an embedding $N\to S^5$ if and only if 
$[f|_{\partial B^3}]=0
\in H_2(S^5-f(\Int N_0))\cong H_1(N_0,\partial N_0)\cong H_1(N)$.}
We conjecture that the Alexander dual of $[f|_{\partial B^3}]$ equals to 
$\pm BH(f')$, where $f':N\to\R^6$ is any extension of $f$.

For an embedding $f:N\to\R^5$ denote by $w(f)\in H_2(N;\Z_2)$ be the Whitney 
invariant [Sk07, \S2]; $s_f\in H_2(N;\Z_2)$ is defined in \S1. 
It would be interesting to know whether $w(f)=s_f$ (this equality would imply 
a simple direct proof of the Relation Lemma from \S3 because $W(if)=\beta w(f)$ 
[Sk05, the Whitney Invariant Lemma ($\beta$)]). 
The equality $w(f)=s_f$ is clear when $f$ is regular homotopic to $g$.  
Indeed, then the regular homotopy together with Seifert surfaces of $f$ and 
$g$ form an immersion $F:V\to\R^6$ of a 4-manifold $V$ such that 
$w(f)=[\Sigma(F)]=w_2(\nu_F)=s_f$. 
For the general case the Sz\"ucs formula for $[\Sigma(F)]+s_f$ in terms of 
singular points of $F$ could be useful. 

%[Hi61, Ki89].
%{\it Proof of the `if' in the Compression Remark.}
%Suppose that $u\in H_1(N)$ and $2u=0$.
%The group $H_1(N)$ is in 1--1 correspondence with the set of regular homotopy
%classes of immersions $N_0\to S^5$ (analogously to [Vra89, p. 167] using
%[Hir59]).
%Since $N_0$ has a 2-dimensional spine, by general position an immersion
%corresponding to $u$ is regularly homotopic to an embedding $f:N_0\to S^5$.
%Then $[f|_{\partial B^3}]=2u=0$ by [Vr89, Addendum 2.2].
%Hence by the above conjecture $f$ extends to an embedding $\bar f:N\to S^5$.
%We have $W(\bar f)=u$ by [Vr89, Theorem 3.1].
%The required results of [Vr89] are proved using [Vr89, \S1, second assertion of Lemma 2.11,
%Lemma 2.12, Proposition 2.8] and so do not use the assumption $k\ge1$.

The higher-dimensional PL analogue of the Compression Theorem (1) is as 
follows [Vr89, Corollary 3.2].

{\it Let $N$ be an $n$-dimensional closed $k$-connected manifold PL embeddable
into $\R^{2n-k-1}$ (the PL embeddability is equivalent to the triviality of the 
normal Stiefel-Whitney class $\overline W_{n-k-1}$).
Suppose that $1\le k\le n-4$.
Then the image of the composition
$\Emb_{PL}^{2n-k-1}(N)\to\Emb_{PL}^{2n-k}(N)\overset{W_{PL}}\to\to H_{k+1}(N)$
is the subgroup formed by elements of order 2 in $H_{k+1}(N)$.
Here the coefficients are $\Z$ for $n-k$ odd and $\Z_2$ for $n-k$ even, and
$W_{PL}$ is bijective.}

%(which was wrongly called the 2-torsion subgroup in [Vr89, Corollary 3.2]).
%Then two embeddings $N\to\R^{2n-k-1}$ are isotopic in $S^{2n-k}$ if and only if their
%restrictions to $N-\Int B^n$ are isotopic in $S^{2n-k-1}$, `if' is clear by [Vr77].

The proofs of this result or its smooth analogues does not work for 3-manifolds 
in $\R^5$ because we cannot apply the smooth analogue of 
Penrose-Whitehead-Zeeman Embedding Theorem, and because we do not assume the 
3-manifold $N$ to be simply-connected and so the smooth analogue of 
[Vr89, Theorem 2.1] is false.

Note that the Compression Theorem (2) holds for 3-manifolds $N$
such that every class in $H_2(N)$ is realizable by a disjoint union of embedded 
2-spheres. 
We conjecture that such manifolds are exactly those of the Compression Theorem 
(2). 

\smallskip
{\it Sketch of the Kreck proof that for each spin structure on $N$ there exists 
a framed embedding $N\to S^5$ inducing this spin structure.} 
(Cf. Realization Lemma of \S3.) 
%; by the formula for the Whitney invariant in 
%Remark (9) this implies the Compression Theorem (1).) 
Take given spin sructure on $N$ and the corresponding spin structure on 
$N\times D^2$. 
Since $\sigma/16:\Omega^{spin}_4(S^1)\to\Z$ is an isomorphism   
%[? ne Ki89!] 
and $\sigma(N\times S^1)=p_1(N\times S^1)/3=0$, it follows that there 
exists a spin 5-manifold $X$ with boundary $\partial X=N\times S^1$ (on which 
the restriction spin structure coincides with the prescribed) and 
a map $F:X\to S^1$ extending the projection from the boundary. 
Making spin ($=BO\left<2\right>=BO\left<4\right>$) surgery we may assume that 
$H_2(\partial X)\to H_2(X)$ is an epimorphism, $F_*:H_1(X)\to H_1(S^1)$ is an 
isomorphism and $\pi_1(N\times x)\to\pi_1(X)$ is zero. 
% (or even $\pi_1(W)=0$?) [?]. 
Then using the Mayer-Vietoris sequence and the van Kampen Theorem 
we prove that $\Sigma:=N\times D^2\bigcup\limits_{N\times S^1}X$ 
is a homotopy 5-sphere. 
Hence $\Sigma$ is diffeomorphic to $S^5$ and we are done.

\smallskip
There is a misprint in [Ki89, Corollary 6 in XI.3]: instead of $K$ it should 
be $F$. 

There is a misprint in [Ki89, Corollary 7 in XI.3]: either $F$ should be 
a sphere, or to the right-hand side of the formula the term 
$8\varphi(M,F)$ should be added (indeed, otherwise the formula is wrong for 
$M=Q-\Int B^4$, where $Q$ is a closed 4-manifold).

\bigskip 
{\bf Constructions of invariants.}

\smallskip
(6) Note that in the proof of the second equality of the Addendum 
we essentially proved that 
$$\eta(f\#f_1,f'\#f_1')\equiv\eta(f,f')+\eta(f_1,f_1')\mod d(W(f\#f_1)).$$ 
It would be interesting to see how the Kreck invariants depends on the choice 
of orientations on $N$ and on $\R^6$ 
(and on self-diffeomorphisms $N\to N$ and/or $\R^6\to\R^6$). 

\smallskip
(7) It would be interesting to construct {\it absolute} Kreck invariant 
$\eta(f)$. 
A possible approach to such a construction is as follows. 
Since $H_4(C_f,\partial C_f)\cong[C_f,\C P^\infty]$, there is a connected 
oriented proper 4-submanifold $X\subset C_f$ with boundary such that 
$[X]=2A_f$. 
We can set $\eta_X(f):\equiv\sigma(X)/16\mod d(W(f))$. 
This $\eta_X(f)$ need not be an integer. 
%It would be interesting to know whether the fractional part of $\eta(f)$ could 
%ever be non-zero and if yes, which meaning it has.)    
The residue $\eta_X(f)$ is independent of change of $X$ by adding a boundary 
(note that analogously to the Twisting Lemma (y) of \S2 $\sigma(X)$ does depend on 
$X$). 
%, so the argument for the independence of $\eta(f)$ on $X$ could not be 
%trivial; 
It would be interesting either to prove that $\eta_X(f)$ is independent of 
$X$, i.e. is independent of adding to $[X]$ a class represented by 
4-submanifold $X_1$ (compact oriented connected, posibly with boundary) of 
$\partial C_f$, or, rather, present 
some additional restrictions on $X$ so that this independence would hold. 
I can only prove that for given $\partial X_1$ and $X$ the residue 
$\eta_{[X]}(f)-\eta_{[X\cup X_1]}(f)\mod16d(W(f))$ is independent of $\Int X_1$, 
but it could probably be non-zero without additional restrictions on $X_1$. 

Compressed Kreck Invariant Lemma (\S3) implies the existence of an {\it absolute 
Kreck invariant} for compressible embeddings $N\to\R^6$.

A related problem is to find for which non-spin $\varphi$ the Kreck invariant 
still equals to $\sigma_{2y}(M_\varphi)/16\mod d(W(f))$ (proof of the 
Codimension Two Compression Theorem suggests that there could be such 
$\varphi$). 

In the construction of $\eta$ instead of fixing an embedding $f'$ we can 
fix a spin simply-connected 4-manifold $W$ such that $\partial W=N$ together 
with an element $y\in H_4(W,\partial W)$ such that 
$[\partial y\times*]=\partial A_f$ on $\partial W\times S^2=\partial C_f$, 
and replace everywhere $C_{f'}$ by $W\times S^2$.  

Note that $16\eta(f,f')=\sigma(X')+\sigma(X'')-X'\cap X''\cap(X'+X'')$ 
if the class $2y$ is represented by the sum of even classes of embedded 
4-manifolds $X'$ and $X''$. 

An alternative proof of the divisibility by 24 in the definition of the Kreck 
invariant is as follows. 
Recall that the group $\Omega^{spin}_6(\C P^\infty)$ is generated by maps
$\C P^3\subset\C P^\infty$ and `the bordism half' of $K\times
S^2\overset{\pr_2}\to\to S^2=\C P^1\subset\C P^\infty$, where $K$
is the Kummer surface.
Hence $\im\omega=\left<(4,1),(-24,0)\right>$.
Now the divisibility by 24 follows because $M$ is spin.

\smallskip
{\it Idea of a possible alternative proof of the Independence Lemma for 
fixed $\varphi$.}
The idea is to use the 4-submanifolds $X\subset M$ and $X_1\subset\partial C_f$ 
representing $2y$ and $2y_1$ instead of $Y$ and $Y_1$. 
Then we try to replace $X\cup X_1$ by an {\it embedded} submanifold in the 
same homology class and calculate its signature. 

Note that in the Cobordism Lemma $M_{\varphi'}$ is not necessarily diffeomorphic  
to $M_\varphi\#\C P^3$. Indeed,  
$M_{\varphi'}\cong(M_\varphi-S^2\times B^3_+\times I)
\bigcup\limits_{\widehat\alpha}S^2\times B^3_-\times I,$ where $\widehat\alpha$ 
is defined analogously to the proof of the Cobordism Lemma.

\smallskip
(8) We conjecture that {\it if $N$ is a closed connected 
orientable $(4k-1)$-manifold, $\overline w_{2k}(N)=0$ and $f:N\to\R^{6k}$ is an 
embedding, then $f\#pt$ is isotopic to $f\#qt$ if and only if $p-q$ is 
divisible by $\frac12d(BH(f))$.}

High connectedness is required in the Higher-dimensional Classification 
Theorem (b) for the above Lemma.  
For $k=2$ we have $\overline w_4(N)=0$ without the connectivity assumption 
[Ma60]. 

The Whitney invariant can analogously be defined for
embeddings of $(2n-m)$-connected closed $n$-manifolds into $\R^m$.
But for $m-n$ even only $W(f)\mod2$ is independent of the isotopy making $f=g$ 
outside $B^n$, because for the equivalent definition in \S1 we 
have $\partial\Sigma(H)=0$ only $\mod2$ (since $H$ is not necessarily an 
immersion).

Analogously the Whitney invariant could be defined for an
embedding $g:N\to\R^6$ such that $g(N)\not\subset\R^5$, but this
is less convenient e.g. because in the main result we would have
$d(u-W(g))$ instead of $d(u)$.

It would be interesting to prove that
$BH(f)-BH(g)=\pm2W(f,g)$ for an embedding $f:N\to\R^{2n-k}$ of a
closed oriented $n$-manifold $N$ and $n-k$ odd. 
Here $W(f,g)$ is defined as in the equivalent definition in \S1. 

It would be interesting to define $BH(F)$ for an isotopy
$F:N^3\times I\to\R^6\times I$ and to express $\eta(\overline F)$ in terms of 
$BH(F)$ as a relative analogue of [BH70, Theorem 2.1].

\smallskip
(9) M. Kreck conjectured the following formula for the Whitney invariant. 
(The advantage of this formula is that it does not use an isotopy making
$f(N_0)\subset S^5$, but we would anyway need such an isotopy in other places.)

Fix a stable normal framing $s$ of some embedding $N_0\to\R^8$ coming
from $\R^5$ (this is a stable spin structure on $N$).
A normal framing of $f|_{N_0}$ is called {\it $s$-spin}, if its sum with the
standard framing of $S^6\subset S^8$ is isotopic to $s$.
%(The reason for such a name would be clear from the sequel; 
($s$-spin normal framing of $f|_{N_0}$ should not be confused with a spin structure 
in $\nu(f|_{N_0})$.)

Analogously to \S2, {\it for each stable framing $s$ on $N_0$ and embedding 
$f:N\to\R^6$ there exists a unique (up to isotopy) $s$-spin normal framing of 
$f|_{N_0}$.}

Let $\xi_0:N_0\to\partial C_f\subset C_f$ be the section formed by first 
vectors of the spin framing. 
We conjecture that $W(f)=\xi_0^!\partial A_f$, cf. the Framing Lemma of \S2. 
(If we take as the initial framing $s$ a framing not 
coming from $\R^5$, then we conjecture that $W(f)=\xi_0^!\partial(A_f-A_g)$.)
An equivalent formulation: {\it if a section $\xi$ of $\nu$ extends to a 
framing $\xi$ of $f$ over $N_0$ that is isotopic to a framing 
$\bar g:N_0\times D^3\to\R^6$ of $g|_{N_0}$ such that 
$\bar g(N_0\times D^2)\subset\R^5$, then $\xi^!\partial A_f=\pm W(f)$.}
These conjectures give an equivalent definition of the Whitney invariant. 
If $H_1(N)$ has no 2-torsion, then these conjectures hold for {\it any}
framing by the Framing Lemma. 

{\it Although for the above construction we only need a section of the spin 
framing,
this section $\xi_0$ cannot be defined without defining the spin framing $\xi$.}
Indeed, define $\varphi:S^1\times S^2\times S^2\to S^1\times S^2\times S^2$ by
$(a,b,x)\to (a,b,\varphi(a)x)$, where $\varphi:S^1\to SO_2\to SO_3$ is a
homotopy nontrivial map.
Then $\varphi$ is a diffeomorphism fiberwise over $S^1\times S^2$
(i.e. $\varphi$ defines change of a section in the normal bundle of an
embedding $S^1\times S^2\to\R^6$), preserving
a section but not preserving the standard stable spin structure.

Note that $\xi^!\partial A_f\ne A_{f,\xi}$ for some $\xi$. 
The analogue of this formula is also false e.g. for the standard embedding 
$S^1\times S^1\to\R^3$.  
It would be interesting to know whether $\varphi\partial A_f=\partial A_{f'}$ 
if $\varphi$ preserves the unlinked section (but is not spin).

%\newpage

\enddocument